\documentclass[reqno, 11pt]{amsart}
\usepackage{amsmath}
\usepackage{latexsym}
\usepackage{amsfonts}
\newtheorem{theorem}{Theorem}
\newtheorem{defn}[theorem]{Definition}

\newtheorem{lem}[theorem]{Lemma}

\newtheorem{Pa}{Paper}[section]
\newtheorem{thm}[Pa]{{\bf Theorem}}

\newtheorem{cor}[Pa]{{\bf Corollary}}

\newtheorem{Pb}[Pa]{{\bf Problem}}

\newcommand{\cF}{\mathcal F}

\newcommand{\C}{{\mathbb C}}
\newcommand{\DD}{{\mathbb D}}
\newcommand{\D}{{\mathbb D}}

\newcommand{\T}{{\mathbb T}}
\newcommand{\TT}{{\mathbb T}}

\newcommand{\PP}{{\mathbb P}}
\newcommand{\ZZ}{{\mathbb Z}}
\newcommand{\I}{{\rm I}}
\newcommand{\cE}{{\mathcal E}}
\newcommand{\cS}{{\mathcal S}}

\newcommand{\bU}{\bf U}
\allowdisplaybreaks[1]

\begin{document}

\title[The higher order Carath\'eodory--Julia theorem]{
The higher order Carath\'eodory--Julia theorem and related boundary
interpolation
problems}

\author{Vladimir Bolotnikov}

\address{Department of Mathematics\\
The College of William and Mary \\
Williamsburg, VA 23187-8795\\
USA}

\email{vladi@math.wm.edu}

\author{ Alexander Kheifets}

\address{Department of Mathematics\\
University of Massachusetts Lowell, \\
Lowell, MA 01854, USA\\
USA}

\email{Alexander\_Kheifets@uml.edu}



\begin{abstract}
The higher order analogue of the classical Carath\'eodory-Julia theorem
on boundary angular derivatives has been obtained in \cite{bk3}.
Here we study boundary interpolation problems for Schur class functions
(analytic and bounded by one in the open unit disk) motivated by that
result.
\end{abstract}

\maketitle

\section{Introduction}
\setcounter{equation}{0}

We denote by $\cS$ the Schur class of analytic functions mapping
the open unit disk $\DD$ into its closure. A well known property of Schur
functions $w$ is that the kernel
\begin{equation}
K_w(z,\zeta)=\frac{1-w(z)\overline{w(\zeta)}}{1-z\bar{\zeta}}
\label{1.1b}
\end{equation}
is positive on $\DD\times\DD$ and therefore, that the matrix
\begin{equation}
{\bf P}^w_{n}(z):=\left[\frac{1}{i!j!} \,
\frac{\partial^{i+j}}{\partial z^{i}\partial\bar{z}^{j}} \,
\frac{1-|w(z)|^2}{1-|z|^2}\right]_{i,j=0}^{n}
\label{1.1}
\end{equation}
which will be referred to as to a {\em Schwarz-Pick matrix},
is positive semidefinite for every $n\ge 0$ and $z\in\D$. We extend  this
notion to boundary points as follows: {\em given a
point $t_0\in\T$, the boundary Schwarz-Pick matrix is
\begin{equation}
{\bf P}^w_{n}(t_0)=\lim_{z\to t_0}{\bf P}^w_{n}(z),
\label{1.2}
\end{equation}
provided the limit in $(\ref{1.2})$ exists}. It is clear that
once the boundary Schwarz-Pick matrix ${\bf P}^w_{n}(t_0)$ exists
for $w\in\cS$, it is positive semidefinite. In \eqref{1.2} and in what
follows, all the limits are nontangential, i.e., $z\in\DD$ tends to a
boundary point nontangentially.
Let us assume that $w\in\cS$ possesses nontangential boundary limits
\begin{equation}
w_j(t_0):=\lim_{z\to t_0}\frac{w^{(j)}(z)}{j!}\quad\mbox{for} \; \;
j=0,\ldots,2n+1
\label{1.3}
\end{equation}
and let
\begin{equation}
\PP^w_{n}(t_0):=\left[\begin{array}{ccc} w_1(t_0) & \cdots &
w_{n+1}(t_0) \\ \vdots & &\vdots \\
w_{n+1}(t_0) & \cdots & w_{2n+1}(t_0)\end{array}\right]{\bf \Psi}_n(t_0)
\left[\begin{array}{ccc}{w}_0(t_0)^* & \ldots &
{w}_{n}(t_0)^*\\ & \ddots & \vdots \\ 0 &&{w_0}(t_0)^*
\end{array}\right],
\label{1.5}
\end{equation}
where the first factor is a Hankel matrix, the third factor
is an upper triangular Toeplitz matrix and where ${\bf
\Psi}_n(t_0)=\left[\Psi_{j\ell}\right]_{j,\ell=0}^n$ is the upper
triangular matrix
\begin{equation}
{\bf \Psi}_n(t_0)=\left[\begin{array}{ccccc} t_0 & -t_0^2& t_0^3&
\cdots& (-1)^{n}{\scriptsize\left(\begin{array}{c} n \\ 0
\end{array}\right)}t_0^{n+1}
\\ 0 & -t_0^3 & 2t_0^4& \cdots &
(-1)^{n}{\scriptsize\left(\begin{array}{c} n \\ 1
\end{array}\right)}t_0^{n+2}\\ \vdots&& t_0^5 & \cdots &
(-1)^{n}{\scriptsize\left(\begin{array}{c}
n \\ 2 \end{array}\right)}t_0^{n+3}\\
\vdots& & &\ddots &\vdots
\\ 0 & \cdots &\cdots & 0 & (-1)^{n}{\scriptsize\left(\begin{array}{c}
n \\ n \end{array}\right)}t_0^{2n+1}\end{array}\right],
\label{1.4}
\end{equation}
with entries
\begin{equation}
\Psi_{j\ell}=\left\{\begin{array}{ccl}
0, & \mbox{if} & j>\ell,
\\ (-1)^{\ell}{\scriptsize
\left(\begin{array}{c} \ell \\ j
\end{array}\right)}t_0^{\ell+j+1}, & \mbox{if} & j\leq\ell.
\end{array}\right.
\label{1.6}
\end{equation}
For notational convenience, in \eqref{1.5} and in what follows we use the
symbol $a^*$ for the complex conjugate of $a\in\C$.

We denote the lower diagonal entry in the
Schwarz-Pick matrix ${\bf  P}^w_{n}(z)$ by
\begin{equation}
d_{w,n}(z):=\frac{1}{(n!)^2}\frac{\partial^{2n}}{\partial
z^{n}\partial\bar{z}^{n}} \, \frac{1-|w(z)|^2}{1-|z|^2}.
\label{1.7}
\end{equation}
The following theorem was obtained in \cite{bk3}.
\begin{thm}
For $w\in\cS$, $t_0\in\TT$ and $n\in{\mathbb Z}_+$, the following are
equivalent:
\begin{enumerate}
\item The following limit inferior is finite
\begin{equation}
\liminf_{z\to t_0}d_{w,n}(z)<\infty
\label{0.1}
\end{equation}
where $z\in\DD$ approaches $t_0$ unrestrictedly.
\item The following nontangential boundary limit exists and is finite:
\begin{equation}
d_{w,n}(t_0):={\displaystyle \lim_{z\to
t_0}d_{w,n}}(z)<\infty.\label{0.18m}
\end{equation}
\item The boundary Schwarz-Pick matrix ${\bf P}^w_{n}(t_0)$
defined via the nontangential boundary limit \eqref{1.2} exists.
\item The nontangential boundary limits $(\ref{1.3})$ exist and satisfy
\begin{equation}
|w_0(t_0)|=1\quad\mbox{and}\quad \PP^w_{n}(t_0)\ge
0,\label{0.13n}
\end{equation}
where $\PP^w_{n}(t_0)$ is the matrix defined in
$(\ref{1.5})$.
\end{enumerate}
Moreover, when these conditions hold, then
\begin{equation}
{\bf P}^w_{n}(t_0)=\PP^w_{n}(t_0).
\label{0.15a}
\end{equation}
\label{T:1.2aa}
\end{thm}
In case $n=0$, Theorem \ref{T:1.2aa} reduces to the classical
Carath\'eodory-Julia theorem \cite{Carat,julia}; this has been discussed
in detail in \cite{bk3}. The relation
$$
d_{w,n}(t_0)=\left[\begin{array}{ccc}w_{n+1}(t_0) & \cdots &
w_{2n+1}(t_0)\end{array}\right]{\bf \Psi}_n(t_0)\left[\begin{array}{c}
{w}_{n}(t_0)^* \\ \vdots \\
{w}_{0}(t_0)^*\end{array}\right]
$$
expresses equality of the lower diagonal entries in \eqref{0.15a};
upon separating the term corresponding to $i=n$ it can be written as
\begin{equation}
d_{w,n}(t_0)=\sum_{i=0}^{n-1}\sum_{j=0}^n
w_{n+i+1}(t_0)\Psi_{ij}(t_0){w}_{n-j}(t_0)^*
+(-1)^nt_0^{2n+1}w_{2n+1}(t_0){w}_0(t_0)^*.
\label{0.0}
\end{equation}
Theorem \ref{T:1.2aa} motivates the following interpolation problem:
\begin{Pb}
Given points $t_1,\ldots,t_k\in\TT$, given integers $n_1,\ldots,n_k\ge
0$ and given numbers $c_{i,j}$ ($j=0,\ldots,2n_i+1; \; i=1,\ldots,k)$,
find all Schur functions $w$ such that
\begin{equation}
\liminf_{z\to t_i}d_{w,n_i}(z)<\infty\quad (i=1,\ldots,k)
\label{up}
\end{equation}
and
\begin{equation}
w_j(t_i):=\lim_{z\to t_i}\frac{w^{(j)}(z)}{j!}=c_{i,j}\quad (i=1,\ldots,k;
\;
j=0,\ldots,2n_i+1).
\label{8.8a}
\end{equation}
\label{P:8.1}
\end{Pb}
The problem makes sense since conditions (\ref{up})
guarantee the existence of the nontangential limits (\ref{8.8a});
upon  preassigning the values $w_j(t_i)$ for $i=1,\ldots,k$ and
$j=0,\ldots,2n_i+1$ , we come up with interpolation Problem \ref{P:8.1}.
It is convenient to reformulate  Problem \ref{P:8.1} in the following form:
\begin{Pb}
Given points $t_1,\ldots,t_k\in\TT$, given integers $n_1,\ldots,n_k\ge
0$ and given numbers
$$
c_{i,j}\quad\mbox{and}\quad\gamma_i \quad (j=0,\ldots,2n_i; \;
i=1,\ldots,k),
$$
find all Schur functions $w$ such that
\begin{equation}
d_{w,n_i}(t_i):=\frac{1}{(n_i!)^2}\lim_{z\to
t_i}\frac{\partial^{2n_i}}{\partial z^{n_i}\partial\bar{z}^{n_i}}
\, \frac{1-|w(z)|^2}{1-|z|^2}=\gamma_i \label{8.7p}
\end{equation}
and
\begin{equation}
w_j(t_i)=c_{i,j}\quad (i=1,\ldots,k; \; j=0,\ldots,2n_i).
\label{8.8b}
\end{equation}
\label{P:8.2}
\end{Pb}
If $w$ is a solution to Problem \ref{P:8.1}, then
conditions (\ref{up}) guarantee the existence of the
nontangential limits \eqref{8.7p} and by a virtue of \eqref{0.0},
\begin{eqnarray}
d_{w,n_i}(t_i)&=&\sum_{\ell=0}^{n_i-1}\sum_{j=0}^{n_i}
w_{n_i+\ell+1}(t_i)\Psi_{\ell j}(t_0){w}_{n_i-j}(t_i)^*\nonumber\\
&&+(-1)^{n_i}t_i^{2n_i+1}w_{2n_i+1}(t_i){w}_0(t_i)^*.
\label{upa}
\end{eqnarray}
Thus, for every Schur function $w$, satisfying (\ref{up}) and
(\ref{8.8a}), conditions \eqref{8.7p} hold with
\begin{equation}
\gamma_i=\sum_{\ell=0}^{n-1}\sum_{j=0}^{n_i}
c_{i,n_i+\ell+1}\Psi_{\ell j}(t_0){c}_{i,n_i-j}^*
+(-1)^{n_i}t_i^{2n_i+1}c_{i,2n_i+1}{c}_{i,0}^*.
\label{upb}
\end{equation}
Conversely, if $w$ is a solution of Problem \ref{P:8.2},
then it  clearly satisfies \eqref{up} and by Theorem \ref{T:1.2aa}, all
the limits in \eqref{8.8a} exist and satisfy relation \eqref{upa}. Since
$w_0(t_i)$ is {\em unimodular}, the equation (\ref{upa}) can be solved for
$w_{2n_i+1}(t_i)$; on account of interpolation conditions \eqref{8.8b}, we
have
\begin{equation}
w_{2n_i+1}(t_i)=(-1)^n\overline{t}_i^{\ 2n_i+1}\left(d_{w,n_i}(t_i)-
\sum_{\ell=0}^{n_i-1}\sum_{j=0}^{n_i}
c_{i,n_i+\ell+1}\Psi_{\ell j}(t_i){c}_{i,n_i-j}^*\right)c_{i,0}.
\label{8.8e}
\end{equation}
It is readily seen now that $w$ is a solution of Problem \ref{P:8.1} with
the data $c_{i,2n_i+1}$ chosen by
\begin{equation}
c_{i,2n_i+1}=(-1)^n\overline{t}_i^{\ 2n_i+1}\left(\gamma_i-
\sum_{\ell=0}^{n_i-1}\sum_{j=0}^{n_i}
c_{i,n_i+\ell+1}\Psi_{\ell j}(t_i){c}_{i,n_i-j}^*\right)c_{i,0}.
\label{8.8f}
\end{equation}
It is known that boundary interpolation problems become more
tractable if they involve  inequalities. Such a relaxed problem is
formulated below; besides of certain independent interest it will serve as
an important intermediate step in solving Problem \ref{P:8.1}.
\begin{Pb}
Given points $t_1,\ldots,t_k\in\TT$, given integers
$n_1,\ldots,n_k\ge 0$ and given numbers $c_{i,j}$ and $\gamma_i\
(j=0,\ldots,2n_i; \; i=1,\ldots,k)$, find all Schur functions $w$
such that
\begin{equation}
d_{w,n_i}(t_i)\le\gamma_i,\label{8.8c}
\end{equation}
\begin{equation}
w_j(t_i)=c_{i,j}\quad (i=1,\ldots,k; \;
j=0,\ldots,2n_i). \label{8.8cc}
\end{equation}
\label{P:8.3}
\end{Pb}
By Theorem \ref{T:1.2aa}, for every solution $w$ of Problem
\ref{P:8.2} there exists the limit
$w_{2n_i+1}(t_i):={\displaystyle\lim_{z\to
t_i}\frac{w^{(2n_i+1)}(z)}{(2n_i+1)!}}$ which satisfies \eqref{8.8e}.
Let $c_{i,2n_i+1}$ be defined as in (\ref{8.8f}). Then it follows from
 \eqref{8.8e}, (\ref{8.8f}) and \eqref{8.8c} that
\begin{equation}
0\le\gamma_i-d_{w,n_i}(t_i)=(-1)^{n_i}t_i^{2n_i+1}\left(c_{i,2n_i+1}-
w_{2n_i+1}(t_i)\right){c}_{i,0}^*, \label{8.7w}
\end{equation}
It is convenient to reformulate Problem \ref{P:8.3} in the
following equivalent form.
\begin{Pb}
Given the data
\begin{equation}
t_i\in\T\quad\mbox{and}\quad c_{i,j}\in\C
\quad (j=0,\ldots,2n_i+1; \; i=1,\ldots,k),
\label{data}
\end{equation}
find all Schur functions $w$ such that
\begin{equation}
d_{w,n_i}(t_i)\le\gamma_i\ ,\label{8.7}
\end{equation}
\begin{equation}
w_j(t_i)=c_{i,j}\quad (i=1,\ldots,k; \; j=0,\ldots,2n_i) \label{8.8}
\end{equation}
and
\begin{equation}
(-1)^{n_i}t_i^{2n_i+1}\left(c_{i,2n_i+1}-w_{2n_i+1}(t_i)
\right){c}_{i,0}^*\ge 0\quad (i=1,\ldots,k).
\label{8.9}
\end{equation}
where $\gamma_i$'s are defined by $(\ref{upb})$.
\label{P:6-final}
\end{Pb}
In Section 3 we will construct the Pick matrix $P$ in terms of
the interpolation data \eqref{data} (see formulas
\eqref{3.1}--\eqref{3.2} below). Then we will show that Problem
$\ref{P:6-final}$ has a solution if and only if $|c_{i,0}|=1$ for
$i=1,\ldots,k$ and $P\ge 0$. In case $P$ is singular, Problem
$\ref{P:6-final}$ has a unique solution $w$ which is a finite Blaschke
product of degree $r\le {\rm rank} \, P$. This unique $w$ may or may not
be a
solution of  Problem \ref{P:8.1}. The case when $P$ is positive definite
is more interesting.
\begin{thm}
Let $|c_{i,0}|=1$ for $i=1,\ldots,k$ and $P>0$. Then
\begin{enumerate}
\item Problem~$\ref{P:6-final}$ has infinitely many solutions which are
parametrized by the linear fractional transformation
\begin{equation}
w(z)=s_0(z)+s_2(z)\left(1-\cE(z)s(z)\right)^{-1}\cE(z)s_1(z)
\label{1.4b}
\end{equation}
where $\cE$ is a free parameter running over the Schur class
${\mathcal S}$ and where the coefficient matrix
\begin{equation}
{\bf S}(z)=\left[\begin{array}{cc}s_0(z) & s_2(z) \\
s_1(z) & s(z)\end{array}\right]
\label{1.4c}
\end{equation}
is rational and inner in $\DD$.
\item A function $w$ of the form $(\ref{1.4b})$ is a solution of
Problem~$\ref{P:8.1}$ if and only if either
\begin{equation}
\liminf_{z\to t_i} \frac{1-|\cE(z)|^2}{1-|z|^2}=\infty
\quad\mbox{or}\quad \lim_{z\to t_i}\cE(z)\neq  {s(t_i)}^*
\label{1.4cd}
\end{equation}
for $i=1,\ldots,k$, where the latter limit is understood as
nontangential, and $s$ is
the right bottom entry of the coefficient matrix ${\bf S}(z)$.
\end{enumerate}
\label{T:1.5}
\end{thm}
Boundary interpolation problems for Schur class functions closely
related to Problem \ref{P:6-final} were studied previously in
\cite{ball8}--\cite{boldym1}, \cite{Kovalishinabound2}. Interpolation
conditions (\ref{8.8} and (\ref{8.9}) in there were accompanied by
various additional restrictions that in fact are equivalent
to our conditions \eqref{8.7}. Establishing these equivalences is a
special issue which will be presented elsewhere. A version of Problem
\ref{P:8.1} (with certain assumptions on the data that guarantee
(\ref{up}) to be in force) was studied in \cite{bgr} for rational matrix
valued Schur functions. In this case, the parameters $\cE$ in the
parametrization formula (\ref{1.4b}) are also rational and therefore, the
situation expressed by the first relation in (\ref{1.4cd}) does not come
into play. A similar matrix valued problem was considered in
\cite{boldym1} where the solvability criteria was established rather than
the description of all solutions. Problem \ref{P:8.2} was considered in
\cite{Sarasonnp} in case $n_1=\ldots=n_k=0$; the second part in Theorem
\ref{T:1.5} can be considered as a higher order generalization of some
results in \cite{Sarasonnp}.

The paper is organized as follows. In Section 2 we recall some
needed results from \cite{bk3} and present some consequences of
conditions (\ref{8.7}) holding for a Schur class function. In
Section 3 we introduce the Pick matrix $P$ in terms of the
interpolation data and establish the Stein equality this matrix
satisfies. In Section 4 we imbed Problem~\ref{P:6-final} in the
general scheme of the Abstract Interpolation Problem (AIP)
developed in \cite{kky,kheyu,kh}. In Section 5 we recall some
needed results on AIP and then prove the first part of Theorem
\ref{T:1.5} in Section 6. Explicit formulas for the coefficients
in the parametrization formula \eqref{1.4b} are derived in Theorem
\ref{T:H2}. Explicit formula for the unique solution of Problem
\ref{P:6-final} in case $P$ is singular is given in Theorem
\ref{T:H1}. In Section 6 we also prove certain properties of the
coefficient matrix \eqref{1.4c} which enable us to prove the
second part of Theorem \ref{T:1.5} in Section 7.

\section{Preliminaries}
\setcounter{equation}{0}

The proof of Theorem \ref{T:1.2aa} presented in \cite{bk3} relies on
the de Branges-Rovnyak spaces $L^w$ and $H^w$ associated to a Schur
function $w$. In this section we recall some needed definitions
and results.  We use the standard notation  $L_2$ for the Lebesgue space
of square integrable functions on the unit circle $\T$; the symbols
$H_2^+$ and $H_2^-$ stand for the Hardy spaces of functions with vanishing
negative (respectively, nonnegative) Fourier coefficients.
The elements in $H_2^+$ and $H_2^-$ will be identified with their unique
analytic (resp., conjugate-analytic) continuations inside the unit disk
and consequently $H_2^+$ and $H_2^-$ will be identified the Hardy spaces
of the unit disk.

Let $w$ be a Schur function. The nontangential boundary limits $w(t)$
exist and are bounded by one at a.e. $t\in\T$ and the matrix valued
function
${\scriptsize\left[\begin{array}{cc} 1& w(t) \\
{w}(t)^* & 1\end{array}\right]}$ is defined and positive semidefinite
almost everywhere on $\T$. The space $L^w$ is the the range
space ${\scriptsize\left[\begin{array}{cc} 1& w \\
{w}^* & 1\end{array}\right]}^{1/2}(L_2\oplus L_2)$ endowed with
the range norm. The set of functions ${\scriptsize\left[\begin{array}{cc}
1& w \\ {w}^* & 1\end{array}\right]}f$ where $f\in L_2\oplus L_2$
is  dense in  $L^w$ and
\begin{equation}
\left\|\left[\begin{array}{cc}
1& w \\ {w}^* & 1\end{array}\right]f\right\|^2_{L^w}=
\left\langle \left[\begin{array}{cc} 1& w \\ {w}^* &
1\end{array}\right]f, \; f\right\rangle_{L_2\oplus L_2}.
\label{2.1}
\end{equation}
\begin{defn}
A function $f={\scriptsize\left[\begin{array}{c}f_+ \\
f_-\end{array}\right]}$ is said to belong to the de
Branges--Rovnyak space $H^w$ if it belongs to $L^w$ and if
$f_+\in H_2^+$ and $f_-\in H_2^-$.
\label{D:2.1}
\end{defn}
As it was shown in \cite{bk3}, the vector valued functions
\begin{equation}
K_{z}^{(j)}(t)=\frac{1}{j!} \,
\frac{\partial^j}{\partial\bar{z}^j} \left(
\left[\begin{array}{cc} 1& w(t) \\
{w}(t)^* & 1\end{array}\right]\left[\begin{array}{c} 1 \\
-{w}(z)^* \end{array}\right]\ \cdot\ \frac{1}{1-t\bar{z}}\right)
\label{2.2}
\end{equation}
defined for $z\in\DD$, $t\in\TT$ and $j\in{\mathbb Z}_+$, belong to the
space $H^w$ and furthermore, for every $z\in\DD$ and every
$f=\left[\begin{array}{c}f_+ \\ f_-\end{array}\right]\in H^w$,
\begin{equation}
\left\langle f, \; K_{z}^{(j)}\right\rangle_{H^w}=
\frac{1}{j!}\frac{\partial^j}{\partial
z^j}f_+(z).\label{2.3}
\end{equation}
Setting $f=K_{\zeta}^{(i)}$ in \eqref{2.3}, we get
\begin{equation}
\left\langle K_{\zeta}^{(i)}, \;
K_{z}^{(j)}\right\rangle_{H^w}= \frac{1}{j!i!} \,
\frac{\partial^{j+i}}{\partial z^j\partial\bar{\zeta}^i}
\left(\frac{1-w(z)\overline{w(\zeta)}}{1-z\bar{\zeta}}\right).\label{2.4}
\end{equation}
Upon differentiating in (\ref{2.2}) and taking into
account that $|t|=1$, we come to the following explicit formulas
for $K_{z}^{(j)}$:
\begin{equation}
K_{z}^{(j)}(t)=\left[\begin{array}{cc} 1& w(t) \\
{w}(t)^* & 1\end{array}\right]\left[\begin{array}{c}
t^j(1-t\bar{z})^{-j-1} \\
-{\displaystyle\sum_{\ell=0}^j{w}_\ell(z)^*t^{j-\ell}
(1-t\bar{z})^{\ell-j-1}}\end{array}\right].
\label{2.5}
\end{equation}
where $w_\ell(z)$ are the Taylor coefficients from the expansion
$$
w(\zeta)=\sum_{\ell=0}^\infty w_\ell(z) (\zeta-z)^\ell, \quad
w_\ell(z)=\frac{w^{(\ell)}(z)}{\ell!}.
$$
The two next theorems (also proved in \cite{bk3}) explain the role of the
condition \eqref{0.1}.
\begin{thm}
Let $w\in{\mathcal S}$, $t_0\in\TT$, $n\in\ZZ_+$ and let
\begin{equation}
\liminf_{z\to t_0}d_{w,n}(z)<\infty
\label{2.6}
\end{equation}
Then the nontangential boundary limits
\begin{equation}
w_j(t_0):=\lim_{z\to t_0}\frac{w^{(j)}(z)}{j!}\quad\mbox{exist for} \; \;
j=0,\ldots,n \label{2.7}
\end{equation}
and the functions
\begin{equation}
K_{t_0}^{(j)}(t)=\left[\begin{array}{cc} 1& w(t) \\
{w}(t)^* & 1\end{array}\right]\left[\begin{array}{c}
t^j(1-t\bar{t}_0)^{-j-1} \\
-{\displaystyle\sum_{\ell=0}^j{w}_\ell(t_0)^*t^{j-\ell}
(1-t\bar{t}_0)^{\ell-j-1}}\end{array}\right]\label{2.8}
\end{equation}
belong to the space $H^w$ for $j=0,\ldots,n$. Moreover,
the kernels $K_{z}^{(j)}$ defined in $(\ref{2.5})$
converge to $K_{t_0}^{(j)}$ for $j=1,\ldots,n$
in norm of $H^w$ as $z\in\DD$ approaches $t_0$ nontangentially:
$$
K_{z}^{(j)}\stackrel{\ \ H^w}{\longrightarrow}K_{t_0}^{(j)}\quad
\mbox{for}\quad j=1,\ldots,n \quad\mbox{as} \; \; (z\to t_0).
$$
\label{T:3.1}
\end{thm}
The next theorem is sort of converse to Theorem \ref{T:3.1}.
\begin{thm}
Let $w\in{\mathcal S}$, $t_0\in\TT$, $n\in\ZZ_+$.
If the numbers $c_0,\ldots,c_n$ are such that the function
$$
F(t)=\left[\begin{array}{cc} 1& w(t) \\
{w}(t)^* & 1\end{array}\right]\left[\begin{array}{c}
t^n(1-t\bar{t}_0)^{-n-1} \\
-{\displaystyle\sum_{\ell=0}^n{c}_\ell^* t^{n-\ell}
(1-t\bar{t}_0)^{\ell-n-1}}\end{array}\right]
$$
belongs to  $H^w$, then  condition \eqref{2.6} holds, the limits
\eqref{2.7} exist and $w_j(t_0)=c_j$ for $j=0,\ldots,n$;
consequently, $F$ coincides with $K_{t_0}^{(n)}$.
\label{T:3.2}
\end{thm}
Now the preceding analysis can be easily extended to a multi-point
setting. Given a Schur function $w$ and $k$-tuples ${\bf
z}=(z_1,\ldots,z_k)$ of points in $\DD$ and ${\bf n}=(n_1,\ldots,n_k)$
of nonnegative integers, define the  {\em generalized Schwarz-Pick
matrix}
\begin{equation}
{\bf P}^w_{\bf n}({\bf z}):=\left[\left[\frac{1}{\ell!r!} \, \left.
\frac{\partial^{\ell+r}}{\partial z^{\ell}\partial\bar{\zeta}^{r}} \,
\left(\frac{1-w(z)\overline{w(\zeta)}}{1-z\bar{\zeta}}\right)
\right\vert_{{\scriptsize\begin{array}{l}z=z_i,\\
\zeta=z_j\end{array}}}\right]_{{\scriptsize\begin{array}{l}
\ell=0,\ldots,n_i\\ r=0,\ldots,n_j\end{array}}}
\right]_{i,j=1}^k.
\label{2.9}
\end{equation}
Given a tuple ${\bf t}=(t_1,\ldots,t_k)$ of distinct points
$t_i\in\T$, define the boundary generalized Schwarz-Pick matrix
\begin{equation}
{\bf P}^w_{\bf n}({\bf t}):=\lim_{{\bf z}\to {\bf t}}{\bf P}^w_{\bf
n}({\bf z})
\label{2.10}
\end{equation}
provided the latter limit exists, where ${\bf z}\to {\bf t}$ means that
$z_i\in\DD $ approaches $t_i$ for $i=1,\ldots,k$ nontangentially.
It is readily seen that conditions
\begin{equation}
\liminf_{z\to t_i}d_{w,n_i}(z)<\infty\quad\mbox{for} \; \; i=1,\ldots,k,
\label{2.11}
\end{equation}
(where  $d_{w,n_i}$ is defined via formula \eqref{1.7}) are necessary for
the limit \eqref{2.10} to exist. They are also sufficient as the next
theorem shows.
\begin{thm}
Let ${\bf t}=(t_1,\ldots,t_k)$ be a tuple of distinct points $t_i\in\T$,
let ${\bf n}=(n_1,\ldots,n_k)\in{\mathbb Z}_+^k$ and let $w$ be a Schur
function satisfying conditions \eqref{2.11}. Then
\begin{enumerate}
\item The following nontangential boundary limits exist:
\begin{equation}
w_j(t_i):=\lim_{z\to t_i}\frac{w^{(j)}(z)}{j!}\quad
(j=0,\ldots,2n_i+1; \; \; i=1,\ldots,k).
\label{2.12}
\end{equation}
\item The functions
\begin{equation}
K_{t_i}^{(j)}(t)=\left[\begin{array}{cc} 1& w(t) \\
{w}(t)^* & 1\end{array}\right]\left[\begin{array}{c}
t^j(1-t\bar{t}_i)^{-j-1}\\
-{\displaystyle\sum_{\ell=0}^j{w}_\ell(t_i)^*t^{j-\ell}
(1-t\bar{t}_i)^{\ell-j-1}}\end{array}\right]\label{2.13}
\end{equation}
belong to the space $H^w$ for $j=0,\ldots,n_i$ and $i=1,\ldots,k$.
\item The boundary generalized Schwarz-Pick matrix ${\bf P}^w_{\bf
n}({\bf t})$ defined via the nontangential limit \eqref{2.10} exists
and is equal to the Gram matrix
of the set $\{K_{t_i}^{(j)}: \; j=0,\ldots,n_i; \; i=1,\ldots,k\}$:
\begin{equation}
{\bf P}^w_{\bf n}({\bf t}):=\left[\left[\left\langle K_{t_j}^{(r)}, \;
K_{t_i}^{(\ell)}\right\rangle_{H^w}\right]_{{\scriptsize\begin{array}{l}
\ell=0,\ldots,n_i\\ r=0,\ldots,n_j\end{array}}}\right]_{i,j=1}^k.
\label{2.14}
\end{equation}
\item The matrix ${\bf P}^w_{\bf n}({\bf t})$ can be expressed in terms
of the nontangential limits \eqref{2.12} as follows:
\begin{equation}
{\bf P}^w_{\bf n}({\bf t})=\left[{\bf P}^w_{ij}\right]_{i,j=1}^k
\label{2.15}
\end{equation}
where ${\bf P}^w_{ij}$ is the $(n_i+1)\times(n_j+1)$ matrix defined by
\begin{equation}
{\bf P}^w_{ij}={\bf H}_{ij}{\bf \Psi}_{n_j}(t_j){\bf W}_j^*,
\label{2.16}
\end{equation}
where ${\bf \Psi}_{n_j}(t_j)$ is defined as in $(\ref{1.4})$,
${\bf W}_j$ is the lower triangular Toeplitz matrix given by
\begin{equation}
{\bf W}_j=\left[\begin{array}{cccc}w_{0}(t_j) & 0 & \ldots & 0 \\
w_{1}(t_j)& w_0(t_i) & \ddots & \vdots \\
\vdots& \ddots & \ddots & 0 \\
w_{n_j}(t_j)&  \ldots & w_{1}(t_j) & w_0(t_j)
\end{array}\right]\label{2.17}
\end{equation}
and where ${\bf H}_{ij}$ is the matrix with the entries
\begin{eqnarray}
\left[{\bf H}_{ij}\right]_{r, s}&=&
\sum_{\ell=0}^{r} (-1)^{r-\ell}\left(\begin{array}{c}s+r-\ell \\
s\end{array}\right)\frac{w_{\ell}(t_i)}
{(t_i-t_j)^{s+r-\ell+1}}\nonumber \\
&&-\sum_{\ell=0}^{s} (-1)^{r}\left(\begin{array}{c}s+r-\ell \\
r\end{array}\right)\frac{w_{\ell}(t_j)}{(t_i-t_j)^{s+r-\ell+1}}
\; .\label{2.18}
\end{eqnarray}
if $i\neq j$, and it is the Hankel matrix
\begin{equation}
{\bf H}_{jj}=\left[\begin{array}{cccc} w_{1}(t_j) & w_{2}(t_j) & \cdots &
w_{n_j+1}(t_j)
\\    w_{2}(t_j) & w_{3}(t_j) & \cdots & w_{n+2}(t_j) \\ \vdots & \vdots
&& \vdots\\ w_{n_j+1}(t_j) & w_{n_j+2}(t_i) & \cdots &
w_{2n_j+1}(t_j)\end{array}\right]\label{2.19}
\end{equation}
otherwise.
\end{enumerate}
\label{T:2.4}
\end{thm}
\begin{proof} The two first statements follow by Theorems
\ref{T:1.2aa} and \ref{T:3.1}. Due to relation \eqref{2.4},
the matrix in \eqref{2.9} can be written as
\begin{equation}
{\bf P}^w_{\bf n}({\bf z}):=\left[\left[\left\langle K_{z_j}^{(r)}, \;
K_{z_i}^{(\ell)}\right\rangle_{H^w}\right]_{{\scriptsize\begin{array}{l}
\ell=0,\ldots,n_i\\ r=0,\ldots,n_j\end{array}}}\right]_{i,j=1}^k.
\label{2.20}
\end{equation}
By Statement 3 in Theorem \ref{T:3.1},
$$
K_{z_i}^{(j)}\stackrel{\ \ H^w}{\longrightarrow}K_{t_i}^{(j)}\quad
\mbox{for}\quad j=1,\ldots,n_i; \; \; i=1,\ldots,k
$$
as $z_i$ approaches $t_i$ nontangentially. Passing to the limit in
\eqref{2.20} we get the existence of the boundary generalized
Schwarz-Pick matrix ${\bf P}^w_{\bf n}({\bf t})$ and obtain its
representation \eqref{2.14}. Let us consider the block partitioning
$$
{\bf P}^w_{\bf n}({\bf z})=\left[{\bf P}^w_{ij}(z_i,z_j)\right]_{i,j=1}^k
$$
conformal with that in \eqref{2.15} so that
\begin{equation}
{\bf P}^w_{ij}(z_i,z_j):=\left[\frac{1}{\ell!r!} \, \left.
\frac{\partial^{\ell+r}}{\partial z^{\ell}\partial\bar{\zeta}^{r}} \,
\left(\frac{1-w(z)\overline{w(\zeta)}}{1-z\bar{\zeta}}\right)
\right\vert_{{\scriptsize\begin{array}{l}z=z_i,\\
\zeta=z_j\end{array}}}\right]_{{\scriptsize\begin{array}{l}
\ell=0,\ldots,n_i\\ r=0,\ldots,n_j\end{array}}}.
\label{2.21}
\end{equation}
The direct differentiation in \eqref{2.21} gives
\begin{eqnarray}
\left[{\bf P}^w_{ij}(z_i,z_j)\right]_{\ell,r}&=&
\sum_{s=0}^{{\rm min}\{\ell,r\}}
\frac{(\ell+r-s)!}{(\ell-s)!(r-s)!}
\frac{z_i^{r-s}\bar{z}_j^{\ell-s}}{(1-z_i\bar{z}_j)^{\ell+r-s+1}}\nonumber\\
&-&\sum_{\alpha=0}^\ell\sum_{\beta=0}^r
\sum_{s=0}^{{\rm min}\{\alpha,\beta\}}
\frac{(\alpha+\beta-s)!}{(\alpha-s)!(\beta-s)!}
\frac{z_i^{\beta-s}\bar{z}_j^{\alpha-s}w_{\ell-\alpha}(z_i){w}_{r-\beta}(z_j)^*}
{(1-z_i\bar{z}_j)^{\alpha+\beta-s+1}}.\nonumber
\end{eqnarray}
For $i\neq j$, we pass to the limit in the latter equality as $z_i\to
t_i$ and $z_j\to t_j$ and take into account \eqref{2.12}:
\begin{eqnarray}
\left[{\bf P}^w_{ij}\right]_{\ell,r}&=&
\sum_{s=0}^{{\rm min}\{\ell,r\}}
\frac{(\ell+r-s)!}{(\ell-s)!(r-s)!}
\frac{t_i^{r-s}\bar{t}_j^{\ell-s}}{(1-t_i\bar{t}_j)^{\ell+r-s+1}}\nonumber\\
&-&\sum_{\alpha=0}^\ell\sum_{\beta=0}^r
\sum_{s=0}^{{\rm min}\{\alpha,\beta\}}
\frac{(\alpha+\beta-s)!}{(\alpha-s)!(\beta-s)!}
\frac{t_i^{\beta-s}\bar{t}_j^{\alpha-s}w_{\ell-\alpha}(t_i)
{w}_{r-\beta}(t_j)^*}
{(1-t_i\bar{t}_j)^{\alpha+\beta-s+1}}.\nonumber
\end{eqnarray}
Verification of the fact that the product on the right hand side of
\eqref{2.15} gives the matrix with the same entries, is straightforward
and will be omitted. Finally, it is readily seen from \eqref{2.21}
and \eqref{1.1} that the  $j$-th diagonal block ${\bf P}^w_{jj}(z_j,z_j)$
coincides with the Schwarz-Pick matrix ${\bf P}^w_{n_j}(z_j)$. Therefore,
by Theorem \ref{T:1.2aa} and formula \eqref{1.5}, its nontangential
boundary limit equals
\begin{eqnarray}
{\bf P}^w_{jj}&=&\PP^w_{n_j}(t_j)\label{2.22}\\
&=&\left[\begin{array}{ccc} w_1(t_j) &
\cdots &w_{n_j+1}(t_j) \\ \vdots & &\vdots \\
w_{n_j+1}(t_j) & \cdots & w_{2n_j+1}(t_j)\end{array}\right]{\bf
\Psi}_{n_j}(t_j)\left[\begin{array}{ccc}{w}_0(t_j)^* & \ldots &
{w}_{n_j}(t_j)^*\\ & \ddots & \vdots \\ 0 &&{w}_0(t_j)^*
\end{array}\right],\nonumber
\end{eqnarray}
which coincides with \eqref{2.16} for $j=i$.
\end{proof}

\section{The Pick matrix and the Stein identity}
\setcounter{equation}{0}

The Pick matrix $P$ defined and studied in this section is important for
formulating a solvability criterion for Problem \ref{P:6-final}
and for parametrizing its solution set. The definition of the Pick
matrix is motivated by the formulas for the matrix ${\bf P}_{\bf
n}^w({\bf t})$ discussed in the previous section. Namely,
\begin{equation}
P=\left[P_{ij}\right]_{i,j=1}^k\in\C^{N\times
N}\quad\mbox{where}\quad N=\displaystyle{\sum_{i=1}^k (n_i+1)},
\label{3.1}
\end{equation}
and the block entries $P_{ij}\in\C^{(n_i+1)\times(n_j+1)}$ are
defined by
\begin{equation}
P_{ij} =H_{ij}\cdot {\bf \Psi}_{n_j}(t_j)\cdot W_j^*,\label{3.2}
\end{equation}
where ${\bf \Psi}_{n_j}(t_j)$ is defined as in (\ref{1.4}), where
\begin{eqnarray}
W_i&=& \left[\begin{array}{cccc}c_{i,0} & 0 & \ldots &
0\\ c_{i,1} & c_{i,0} & \ldots & 0\\
\vdots & \ddots & \ddots & \vdots \\
c_{i,n_i} &\ldots &  c_{i,1} & c_{i,0}\end{array}\right],\label{3.3}\\
H_{ii}&=&\left[\begin{array}{cccc} c_{i,1} & c_{i,2} & \cdots &
c_{i,n_i+1}
\\    c_{i,2} & c_{i,3} & \cdots & c_{i, n+2} \\ \vdots & \vdots &&
\vdots\\ c_{i,n_i+1} & c_{i,n_i+2} & \cdots &
c_{i,2n_i+1}\end{array}\right]\label{3.3a}
\end{eqnarray}
for $i=1,\ldots,k$ and where the matrices $H_{ij}$ (for $i\neq j$)
are defined entrywise by
\begin{eqnarray}
\left[H_{ij}\right]_{r, s}&=&
\sum_{\ell=0}^{r} (-1)^{r-\ell}
\left(\begin{array}{c}s+r-\ell \\
s\end{array}\right)\frac{c_{i,\ell}}
{(t_i-t_j)^{s+r-\ell+1}}\nonumber \\
&&-\sum_{\ell=0}^{s} (-1)^{r}\left(\begin{array}{c}s+r-\ell \\
r\end{array}\right)\frac{c_{j,\ell}}{(t_i-t_j)^{s+r-\ell+1}}
\; .\label{3.4}
\end{eqnarray}
for $r=0,\ldots,n_i$ and $s=0,\ldots,n_j$. The latter formulas define
$P$ exclusively in terms of the interpolation data of \eqref{data}. We
also associate with the same data the following matrices:
\begin{align}
&T=\left[\begin{array}{ccc}T_1 & & 0 \\ & \ddots & \\
0 && T_k\end{array}\right],&\quad\mbox{where}\quad &T_i
=\left[\begin{array}{cccc} \bar{t}_i & 1 & \ldots & 0 \\
0 &  \bar{t}_i & & \vdots \\ \vdots & \ddots & \ddots& 1\\
0 & \ldots & 0 & \bar{t}_i\end{array}\right],\label{3.5}\\
&E=\left[\begin{array}{ccc}E_1 & \ldots &
E_k\end{array}\right],&\quad\mbox{where}\quad&
E_i=\left[\begin{array}{cccc} 1& 0& \ldots&
0\end{array}\right],\label{3.6}\\
&M=\left[\begin{array}{ccc}M_1 & \ldots &
M_k\end{array}\right],&\quad\mbox{where}\quad &M_i=
\left[\begin{array}{ccc}{c}_{i,0}^*&
\ldots & {c}_{i,n_i}^*\end{array}\right].\label{3.7}
\end{align}
Note that $T_i\in\C^{(n_i+1)\times (n_i+1)}$ and $E_i, \,
M_i\in\C^{1\times (n_i+1)}$. The main result of
this section is:
\begin{thm}
Let $|c_{i,0}|=1$ for $i=1,\ldots,k$ and let us assume that the
diagonal blocks $P_{ii}$ of the matrix $P$ defined in
$(\ref{3.1})$--$(\ref{3.4})$ are Hermitian for $i=1,\ldots,k$.
Then the matrix $P$ is Hermitian and satisfies the Stein identity
\begin{equation}
P-T^*PT=E^*E-M^*M, \label{3.8}
\end{equation}
where the matrices $T$, $E$ and $M$ are defined in
$(\ref{3.5})$--$(\ref{3.7})$.
\label{T:8.a1}
\end{thm}
In view of (\ref{3.5})--(\ref{3.7}),
verifying (\ref{3.8}) is equivalent to verifying
\begin{equation}
P_{ij}-T_i^*P_{ij}T_j=E_{i}^*E_{j}-M_i^*M_j\quad(i,j=1,\ldots,k).
\label{3.9}
\end{equation}
The identity like that is not totally surprising due to a special
(Hankel and Toeplitz) structure of the factors $H_{ij}$ and $W_j$
in (\ref{3.2}). Indeed, the identity verified in the next lemma
(though, not exactly of the form (\ref{3.9})) follows from the
structure of $P_{ij}$ only (without any symmetry assumptions).
Note that the righthand side in (\ref{3.8}) as well as the one
in (\ref{3.9}) is of rank $2$ .
\begin{lem}
Let $P_{ij}$ be defined as in $(\ref{3.2})$. Then
\begin{equation}
P_{ij}-T_i^*P_{ij}T_j=E_{i}^*\overline{M}_j {\bf
\Psi}_{n_j}(t_j)W_j^*T_j-M_i^*M_j \label{3.10}
\end{equation}
where according to $(\ref{3.7})$,
$$
\overline{M}_j=\left[\begin{array}{cccc} c_{j,0} & c_{j,1} &\ldots &
c_{j,n_j}\end{array}\right]=E_{j}W_j^\top.
$$
\label{L:8.1a}
\end{lem}
\begin{proof} We shall make use of the following equalities
\begin{equation}
W_j^*T_j=T_jW_j^*,\quad \overline{T}_j{\bf
\Psi}_{n_j}(t_j)T_j={\bf \Psi}_{n_j}(t_j),\quad E_{j}{\bf
\Psi}_{n_j}(t_j)T_j=E_{j}. \label{3.11}
\end{equation}
The first equality follows by the Toeplitz triangular structure of $W_j^*$
and $T_j$. The matrix $\overline{T}_j{\bf \Psi}_{n_j}(t_j)T_j$ is
upper triangular as the product of upper triangular matrices and
due to (\ref{1.6}) and (\ref{3.5}), its $s\ell$-th  entry (for
$\ell\ge s$) equals
\begin{eqnarray*}
&&\left[\overline{T}_j{\bf \Psi}_{n_j}(t_j)T_j\right]_{s,\ell}=
\Psi_{s,\ell}+
t_j\Psi_{s,\ell-1}+\bar{t}_j\Psi_{s+1,\ell}+\Psi_{s+1,\ell-1}\\
&&=(-1)^\ell t_j^{s+\ell+1}\left[\left(\begin{array}{c} \ell \\ s
\end{array}\right)-\left(\begin{array}{c} \ell-1 \\ s\end{array}\right)
+\left(\begin{array}{c} \ell \\ s+1\end{array}\right)-
\left(\begin{array}{c} \ell-1 \\ s+1\end{array}\right)\right]\\
&&=(-1)^\ell t_j^{s+\ell+1}\left(\begin{array}{c} \ell \\ s
\end{array}\right)=\Psi_{s,\ell}
\end{eqnarray*}
and completes the verification of the second equality in
(\ref{3.11}). The last relation in (\ref{3.11}) follows by
(\ref{1.6}) and (\ref{3.5}) and (\ref{3.6}):
$$
E_{j}{\bf \Psi}_{n_j}(t_j)T_j=\left[\begin{array}{cccc} t_j &
-t_j^2 & \ldots & (-1)^{n_j}t_j^{n_j+1}\end{array}\right]
T_j=\left[\begin{array}{cccc} 1 & 0 & \ldots &
0\end{array}\right]=E_{j}.
$$
We will also use the identity
\begin{equation}
H_{ij}\overline{T}_j-T_i^*H_{ij}=E_i^*\overline{M}_j-M_i^*E_{j}
\label{3.12}
\end{equation}
which holds for every $i,j=1,\ldots,k$ and is verified by
straightforward calculations (separately for the cases $i=j$ and
$i\neq j$). We have
\begin{eqnarray}
P_{ij}-T_i^*P_{ij}T_j&=& H_{ij}{\bf \Psi}_{n_j}(t_j)W_i^*-
T_i^*H_{ij}{\bf \Psi}_{n_j}(t_j)T_jW_j^*\nonumber\\
&=&\left(H_{ij}\overline{T}_j-T_j^*H_{ij}\right){\bf\Psi}_{n_j}(t_j)
T_jW_j^*\nonumber\\
&=&\left(E_{i}^*\overline{M}_j-M_i^*E_{j}
\right){\bf\Psi}_{n_j}(t_j)T_jW_j^*,\label{3.13}
\end{eqnarray}
where the first equality follows by (\ref{3.2}) and the first
relation in (\ref{3.11}), the second equality relies on the
second relation in (\ref{3.11}) and the last equality is a
consequence of (\ref{3.12}). Combining the third relation in
(\ref{3.11}) with formulas (\ref{3.3}) and (\ref{3.7}) we
get
$$
E_{j}{\bf\Psi}_{n_i}(t_i)T_jW_j^*=E_{j}W_j^*=M_j
$$
which being substituted into (\ref{3.13}) leads us to
(\ref{3.10}).\end{proof}

{\bf Proof of Theorem \ref{T:8.a1}:} By Lemma \ref{L:8.1a} the
structure of $P$ implies (\ref{3.10}). First we consider the case
when $j=i$. Since, by assumption, matrices $P_{ii}$
($i=1,\ldots,k$) are Hermitian, the left hand sides in
(\ref{3.10}) are Hermitian, and hence the right hand sides in
(\ref{3.10}) must be Hermitian. In other words
$$
E_{i}^*\overline{M}_i{\bf \Psi}_{n_i}(t_i)W_i^*T_i=
(\overline{M}_i{\bf \Psi}_{n_i}(t_i)W_i^*T_i)^*E_{i}\quad
\mbox{for} \; \; i=1,\ldots,k.
$$
Multiplying the latter relation by $E_i$ from the left and
taking into account that
$$
E_{i}E_{i}^*=1\quad \mbox{and}\quad E_{i}(\overline{M}_i{\bf
\Psi}_{n_i}(t_i)W_i^*T_i)^*=(c_{i,0}t_ic_{i,0}^*\overline
t_i)^*=1,
$$
we get
$$
\overline{M}_i{\bf \Psi}_{n_i}(t_i)W_i^*T_i= E_{i}\quad
\mbox{for} \; \; i=1,\ldots,k.
$$
Therefore, relations (\ref{3.10}) turn into (\ref{3.9}), which
is equivalent to (\ref{3.8}). Furthermore, for $i\neq j$, the Stein
equation
$$
X-T_i^*XT_j=E_{i}^*E_{j}-M_i^*M_j
$$
has a unique solution $X$. Taking adjoint of both sides in
(\ref{3.9}) we conclude that the matrix $P_{ij}^*$ satisfies the
same Stein equation as $P_{ji}$ does and then, by the above
uniqueness, $P_{ij}^*=P_{ji}$ for $i\neq j$. It follows now that
$P$ is Hermitian. \qed
\begin{thm}
Let $t_1,\ldots,t_k\in\T$, $n_1,\ldots,n_k\in{\mathbb Z}_+$,
$N=\sum_{i=1}^k(n_1+1)$ and let us assume that a Schur function $w$
satisfies conditions \eqref{2.11}. Then the matrix
${\bf P}^w_{\bf n}({\bf t})$ defined via the limit \eqref{2.10} (that
exists by Theorem \ref{3.1}) satisfies the Stein identity
\begin{equation}
{\bf P}^w_{\bf n}({\bf t})-T^*{\bf P}^w_{\bf n}({\bf t})T=E^*E-(M^w)^*M^w,
\label{3.14}
\end{equation}
where the matrices $T$ and $E$ are defined in
$(\ref{3.5})$, $(\ref{3.6})$ and
\begin{equation}
M^w=\left[\begin{array}{ccc}M_1 & \ldots &
M_k\end{array}\right],\quad\mbox{where}\quad M^w_i=
\left[\begin{array}{ccc}{w}_{0}(t_i)^*&
\ldots &{w}_{n_i}(t_i)^*\end{array}\right].\label{3.15}
\end{equation}
\label{T:3.7}
\end{thm}
\begin{proof} By Theorem \ref{T:3.1}, the matrix ${\bf P}^w_{\bf n}({\bf
t})$ admits the representation \eqref{2.15}-\eqref{2.18}, that is the
same
structure as the Pick matrix $P$ constructed in \eqref{3.1}--\eqref{3.4}
but with parameters $c_{ij}$ replaced by $w_j(t_i)$. Furthermore,
it is positive semidefinite (and therefore, its diagonal blocks are
Hermitian) due to representation \eqref{2.14}, whereas $|w_0(t_i)|=1$
for $i=1,\ldots,k$, by Theorem \ref{T:1.2aa}. Upon applying Theorem
\ref{T:8.a1} we conclude that ${\bf P}^w_{\bf n}({\bf t})$
satisfies the same Stein identity as $P$ but with $M^w$ instead of
$M$, i.e., the Stein identity \eqref{3.14}.\end{proof}

\section{Reformulation of Problem \ref{P:6-final}}
\setcounter{equation}{0}

The formula \eqref{2.14} for ${\bf P}_{\bf n}^w({\bf t})$ motivates us
to introduce the matrix function
\begin{equation}
\widetilde{\bf F}^w(t)=\left[\begin{array}{ccc}\widetilde{\bf F}^w_1(t) &
\ldots &\widetilde{\bf F}^w_k(t)\end{array}\right], \label{4.1}
\end{equation}
where
\begin{equation}
\widetilde{\bf F}^w_i(t):=\left[\begin{array}{cccc}K_{t_i}^{(0)}(t) &
K_{t_i}^{(1)}(t) & \ldots &
K_{t_i}^{(n_i)}(t)\end{array}\right]\quad (i=1,\ldots,k),
\label{4.2}
\end{equation}
and $K_{t_i}^{(j)}(t)\ (j=0,\ldots,n_i)$ are the functions
defined in (\ref{2.13}).
\begin{thm}
Let $t_1,\ldots,t_k\in\T$, $n_1,\ldots,n_k\in{\mathbb Z}_+$
and let us assume that a Schur function $w$ satisfies conditions
\eqref{2.11}. Then for $\widetilde{\bf F}^w$ defined in \eqref{4.2},
\eqref{4.3} we have
\begin{enumerate}
\item The function $\widetilde{\bf F}^wx$ belongs to the de
Branges-Rovnyak space $H^w$ for every vector $x\in\C^N$ and
\begin{equation}
\|\widetilde{\bf F}^wx\|_{H^w}^2=x^*{\bf P}_{\bf n}^w({\bf t})x
\label{4.3}
\end{equation}
where ${\bf P}_{\bf n}^w({\bf t})$ is the boundary generalized
Schwarz-Pick matrix (that exists due to conditions \eqref{2.11})
and $N:=\sum_{i=1}^k(n_i+1)$.
\item $\widetilde{\bf F}^w$ admits the representation
\begin{equation}
\widetilde{\bf F}^w(t)=\left[\begin{array}{cc} 1& w(t) \\
{w}(t)^* & 1\end{array}\right]\left[\begin{array}{l}\ \ E \\ -M^w
\end{array}\right]\left(\I-tT\right)^{-1},
\label{4.4}
\end{equation}
\end{enumerate}
where the matrices $T$, $E$ and $M^w$ are defined in \eqref{3.5},
\eqref{3.6} and \eqref{3.15}, respectively.
\label{T:4.1}
\end{thm}
\begin{proof} By Theorem \ref{T:1.2aa},  conditions (\ref{8.7}) guarantee
that the functions $K_{t_i}^{(j)}$ defined in \eqref{2.13} belong to
$H^w$
and the boundary Schwarz-Pick matrix ${\bf P}^w_{\bf n}({\bf t})$ exists
and admits a representation \eqref{2.14}. Now it follows from
\eqref{4.1} and \eqref{4.2} that for every $x\in\C^n$, the function ${\bf
F}^wx$ belongs to $H^w$ as a linear combination of the kernels
$K_{t_i}^{(j)}\in H^w$, while relation \eqref{4.3} is an immediate
consequence of \eqref{2.14}. Furthermore, by definitions (\ref{3.5}),
(\ref{3.6}) and \eqref{3.15} of $T_i$, $E_i$ and $M^w_{i}$,
\begin{equation}
\left[\begin{array}{c} E_{i} \\ -M^w_i
\end{array}\right]\left(\I-tT_i\right)^{-1}=
\left[\begin{array}{ccc}{\displaystyle\frac{1}{1-t\bar{t}_i}} &
\ldots & {\displaystyle\frac{t^{n_i}}{(1-t\bar{t}_i)^{n_i+1}}}\\
-{\displaystyle\frac{{w}_{0}(t_i)^*}{1-t\bar{t}_i}} & \ldots &
-{\displaystyle\sum_{\ell=0}^{n_i}\frac{{w}_{\ell}(t_i)^*t^{n_i-\ell}}
{(1-t\bar{t}_i)^{n_i+1-\ell}}}\end{array}\right].
\label{4.5}
\end{equation}
Multiplying both sides of \eqref{4.5} by the matrix
$\left[\begin{array}{cc}1& w(t) \\ {w}(t)^* & 1\end{array}\right]$
on the left and taking into account (\ref{2.13}) and \eqref{4.2} we get
\begin{eqnarray}
\left[\begin{array}{cc} 1& w(t) \\
{w}(t)^* & 1\end{array}\right]\left[\begin{array}{l}\ \ E \\ -M^w
\end{array}\right]\left(\I-tT\right)^{-1}&=&
\left[\begin{array}{cccc}K_{t_i}^{(0)}(t) &
K_{t_i}^{(1)}(t) & \ldots &
K_{t_i}^{(n_i)}(t)\end{array}\right]\nonumber\\
&=:&\widetilde{\bf F}^w_i(t)\qquad (i=1,\ldots,k).
\label{4.6}
\end{eqnarray}
Now representation formula \eqref{4.4} follows by definitions
(block partitionings) \eqref{4.1}, (\ref{3.5}), (\ref{3.6}) and \eqref{3.15}
of $\widetilde{\bf F}^w$, $T$, $E$ and $M^w$.\end{proof}

Now we modify $\widetilde{\bf F}^w$ replacing $M^w$ by $M$ in \eqref{4.4}:
we introduce the function
\begin{equation}
{\bf F}^w(t):=\left[\begin{array}{cc} 1& w(t) \\
{w}(t)^* & 1\end{array}\right]\left[\begin{array}{l}\ \ E \\ -M
\end{array}\right]\left(\I-tT\right)^{-1}
\label{4.7}
\end{equation}
with $T$, $E$ and $M$ defined in (\ref{3.5})--(\ref{3.7}).
The two next theorems show that Problem \ref{P:6-final} can be
reformulated in terms of this function and of the Pick matrix $P$.
\begin{thm}
Assume that $w$ solves Problem $\ref{P:6-final}$
(i.e., $w\in{\mathcal S}$ and satisfies interpolation conditions
$(\ref{8.7})$--$(\ref{8.9})$) and let ${\bf F}^w$ be defined as in
$(\ref{4.7})$. Then
\begin{enumerate}
\item The function ${\bf F}^w x$ belongs to $H^w$ for every vector
$x\in\C^N$ and
\begin{equation}
\|{\bf F}^w x\|^2_{H^w}\le x^*Px
\label{4.8}
\end{equation}
where $P$ is the Pick matrix defined in  $(\ref{3.1})$--$(\ref{3.4})$.
\item The numbers $c_{i,0}$ are unimodular for $i=1,\ldots,k$ and
the matrix $P$ is positive semidefinite
\begin{equation}
|c_{i,0}|=1 \; \; (i=1,\ldots,k)\quad \mbox{and}\quad P\ge 0.
\label{4.9}
\end{equation}
\item $P$ satisfies the Stein identity \eqref{3.8}.
\end{enumerate}
Furthermore, if $w$ is a solution of  Problem $\ref{P:8.1}$, then
\begin{equation}
\|{\bf F}^w x\|^2_{H^w}= x^*Px\quad\mbox{for every} \; \; x\in\C^N.
\label{4.8a}
\end{equation}

\label{T:4.2}
\end{thm}
\begin{proof} Conditions (\ref{8.7}) guarantee (by Theorem \ref{T:1.2aa})
that the limits $w_0(t_i)$ are unimodular for $i=1,\ldots,k$; since
$w_0(t_i)=c_{i,0}$ (according to \eqref{8.8}), the first condition in
\eqref{4.9} follows.

Conditions (\ref{8.7}) also guarantee (by Theorem \ref{T:4.1}),
that for every $x\in\C^N$, the function $\widetilde{\bf F}^w x$ belongs to
$H^w$ for every vector $x\in\C^N$ and equality \eqref{4.3} holds, where
$\widetilde{\bf F}^w$ is defined by the representation formula \eqref{4.4}.
On account of interpolation conditions \eqref{8.8} (only for $j=0,\ldots,n_i$
and for every $i=1,\ldots,k$) and by definitions \eqref{3.7} and
\eqref{3.15}, it follows that $M=M^w$. Then the formulas \eqref{4.4}
and \eqref{4.7} show that ${\bf F}^w\equiv\widetilde{\bf F}^w$, so that
equality \eqref{4.3} holds with ${\bf F}^w$ instead of $\widetilde{\bf F}^w$:
\begin{equation}
\|{\bf F}^wx\|^2_{H^w}=x^*{\bf P}^w_{\bf n}({\bf t})x.
\label{4.10}
\end{equation}
Thus, to prove \eqref{4.8}, it suffices to show that ${\bf P}^w_{\bf
n}({\bf t})\le P$. We will use formulas \eqref{2.15}--\eqref{2.19}
defining ${\bf P}^w_{\bf n}({\bf t})$ in terms of the boundary limits
$w_j(t_i)$. In view of these formulas and due to interpolation conditions
(\ref{8.8}), ${\bf P}_{\bf n}^w({\bf t})$ can be mostly expressed in terms
of the interpolation data (\ref{data}). Indeed,  comparing
(\ref{3.2})--(\ref{3.4}) and (\ref{2.15})--(\ref{2.18}) we conclude that
\begin{equation}
{\bf P}^w_{ij}=P_{ij}\quad (i\neq j) \label{4.11}
\end{equation}
and that formula (\ref{2.22}) for the diagonal blocks of ${\bf
P}^w$ turns into
\begin{equation}
{\bf P}^w_{ii}=\left[\begin{array}{ccc} c_{i,1} &
\cdots & c_{i,n_i+1}\\    c_{i,2} & \cdots & c_{i,n_i+2}
\\ \vdots && \vdots\\
c_{i,n_i+1} & \cdots & w_{2n_i+1}(t_i)\end{array}\right]{\bf
\Psi}_{n_i}(t_i)\left[\begin{array}{cccc}{c}_{i,0}^* &
{c}_{i,1}^* & \ldots & {c}_{i,n_i}^*\\ 0 & {c}_{i,0}^* & \ldots &
{c}_{i,n_i-1}^* \\ \vdots & \ddots & \ddots & \vdots \\
0 & \ldots & 0 & {c}_{i,0}^*\end{array}\right]. \label{4.12}
\end{equation}
Taking into account the upper triangular structure of ${\bf
\Psi}_{n_i}(t_i)$, we conclude from (\ref{3.2}), (\ref{3.3})
and (\ref{4.12}) that all the corresponding entries in $P_{ii}$
and ${\bf P}^w_{ii}$ are equal except for the rightmost bottom
entries that are equal to  $\gamma_i$ and to
$d_{w,n_i}(t_i)$, respectively. Thus, by condition (\ref{8.7}),
\begin{equation}
P_{ii}-{\bf P}^w_{ii}=\left[\begin{array}{ccc}0 & \ldots & 0 \\
\vdots & \ddots & \vdots \\ 0 & \ldots &
\gamma_i-d_{w,n_i}(t_i)\end{array}\right]\ge 0, \label{4.13}
\end{equation}
for $i=1,\ldots,k$ which together with (\ref{4.11})
imply $P\ge {\bf P}^w$ and therefore, relation \eqref{4.8}.
If $w$ is a solution of Problem  $\ref{P:8.1}$ (or equivalently, of
Problem \ref{P:8.2}), then $\gamma_i-d_{w,n_i}(t_i)=0$ for $i=1,\ldots,k$
in (\ref{4.13}) which proves the final statement in the theorem. Since
${\bf P}^w\ge 0$, we conclude  from the inequality $P\ge {\bf P}^w$ that
$P\ge 0$ which completes the  proof of the second statement of the
theorem. The third statement follows from  \eqref{4.9} by Theorem
\ref{T:8.a1}.
\end{proof}
The next theorem is the converse to Theorem \ref{4.2}.
\begin{thm}
Let $P$, $T$, $E$ and $M$ be the matrices given by
$(\ref{3.1})$--$(\ref{3.7})$. Let $|c_{i,0}|=1$ and $P\ge 0$. Let $w$ be a
Schur function such that
\begin{equation}
{\bf F}^wx:=\left[\begin{array}{cc} 1& w(t) \\
w(t)^* & 1\end{array}\right]\left[\begin{array}{r} E  \\ -M
\end{array}\right]\left(\I-tT\right)^{-1}
x\ \ \mbox{belongs to}\ H^w \label{4.14}
\end{equation}
for every $x\in\C^N$ and satisfies \eqref{4.8}.
Then $w$ is a solution of Problem~$\ref{P:6-final}$. If moreover,
\eqref{4.8a} holds, then $w$ is a solution of Problem~$\ref{P:8.1}$.
\label{T:4.3}
\end{thm}
\begin{proof} By the definitions \eqref{3.5}--\eqref{3.7} of
$T$, $E$ and
$M$, the columns of the $2\times N$ matrix ${\bf F}^w$ defined in
\eqref{4.7}, are of the form
$$
\left[\begin{array}{cc} 1& w(t) \\
{w}(t)^* & 1\end{array}\right]\left[\begin{array}{c}
t^j(1-t\bar{t}_i)^{-j-1} \\
-{\displaystyle\sum_{\ell=0}^j{c}_{i,\ell}^* t^{j-\ell}
(1-t\bar{t}_i)^{\ell-j-1}}\end{array}\right]
$$
for $j=1,\ldots,n_i$ and $i=1,\ldots,k$, and all of them belong to
$H^w$  by the assumption \eqref{4.14} of the theorem. In particular, the
functions
$$
F_i(t)=\left[\begin{array}{cc} 1& w(t) \\
{w}(t)^* & 1\end{array}\right]\left[\begin{array}{c}
t^{n_i}(1-t\bar{t}_i)^{-n_i-1} \\
-{\displaystyle\sum_{\ell=0}^{n_i}{c}_{i,\ell}^* t^{n_i-\ell}
(1-t\bar{t}_i)^{\ell-n_i-1}}\end{array}\right]
$$
belong to $H^w$, which implies, by Theorems \ref{3.2} and \eqref{2.4}
that
\begin{equation}
\liminf_{z\to t_i}d_{w,n_i}(z)<\infty\quad\mbox{for} \; \;
i=1,\ldots,k,
\label{4.15}
\end{equation}
and that the nontangential limits \eqref{2.12} exist and satisfy
\begin{equation}
w_j(t_i)=c_{ij} \quad\mbox{for}\quad j=1,\ldots,n_i \quad\mbox{and}\quad
i=1,\ldots,k.
\label{4.16}
\end{equation}
Therefore, $w$ meets conditions (\ref{8.8}) for
$i=1,\ldots,k$ and $\ell_i=0,\ldots,n_i$. By Theorem \ref{T:4.1},
conditions \eqref{4.15} guarantee that the boundary generalized
Schwarz-Pick matrix  ${\bf P}_{\bf n}^w({\bf t})$ exists and
that
\begin{equation}
\|\widetilde{\bf F}^wx\|_{H^w}^2=x^*{\bf P}_{\bf n}^w({\bf t})x
\quad\mbox{for every} \; \; x\in\C^N,
\label{4.17}
\end{equation}
where $\widetilde{\bf F}^w$ is the $2\times N$
matrix function defined in \eqref{4.4}. By Theorem \ref{T:2.4},
${\bf P}_{\bf n}^w({\bf t})$ is represented in terms of the boundary
limits \eqref{2.12} by formulas (\ref{2.15})--(\ref{2.18}).
Equalities \eqref{4.16} along with definitions \eqref{3.7} and
\eqref{3.15} of $M$ and $M^w$ show that the two latter matrices are
equal and thus ${\bf F}^w\equiv\widetilde{\bf F}^w$, by \eqref{4.4} and
\eqref{4.7}.
Now combining \eqref{4.17} and \eqref{4.14} gives ${\bf P}_{\bf
n}^w({\bf t})\le P$ which implies inequalities for the diagonal blocks
\begin{equation}
{\bf P}_{ii}^w\le P_{ii} \quad (i=1,\ldots,k).
\label{4.18}
\end{equation}
Since $d_{w,n_i}(t_i)$ and $\gamma_i$ are (the lower) diagonal entries
in ${\bf P}_{ii}^w$ and $P_{ii}$ respectively, the latter inequality
implies \eqref{8.9}.

By Theorems \ref{T:8.a1} and \ref{T:3.7}, the matrices $P$ and ${\bf
P}_{\bf n}^w({\bf t})$ possess the Stein identities
\eqref{3.8} and \eqref{3.14}, respectively; since $M=M^w$, the matrix
$\widetilde{P}:= P-{\bf P}_{\bf n}^w({\bf t})$ satisfies the homogeneous
Stein identity
$$
\widetilde{P}-T^*\widetilde{P}T=0.
$$
By the diagonal structure \eqref{3.5} of $T$ and in view of
\eqref{4.18} we have for the diagonal blocks $\widetilde{P}_{ii}$ of
$\widetilde{P}$,
\begin{equation}
\widetilde{P}_{ii}-T_i^*\widetilde{P}_{ii}T_i=0\quad\mbox{and}\quad
\widetilde{P}_{ii}\ge 0 \quad(i=1,\ldots,k).
\label{4.19}
\end{equation}
By the Jordan structure \eqref{3.5} of $T_i$, it follows from
\eqref{4.19}  that $\widetilde{P}_{ii}$ is necessarily of the form
\begin{equation}
\widetilde{P}_{ii}=P_{ii}-{\bf P}^w_{ii}=\left[\begin{array}{ccc}
0 & \ldots & 0 \\
\vdots & \ddots & \vdots \\ 0 & \ldots &
\delta_i\end{array}\right]\quad\mbox{with}\quad \delta_i\ge 0
\label{4.20a}
\end{equation}
(for a simple proof see e.g., \cite[Corollary 10.7]{boldym1}).
On the other hand, by the representations \eqref{2.16} and \eqref{3.2},
$$
{\bf P}^w_{ii}={\bf H}_{ii}{\bf \Psi}_{n_i}(t_i){\bf
W}_i^*\quad\mbox{and}\quad P_{ii}=H_{ii}{\bf \Psi}_{n_i}(t_i)W_i^*
$$
and since by \eqref{4.16}, ${\bf W}_i=W_i$ (which is readily
seen from the definitions \eqref{2.17} and \eqref{3.3}), we conclude
that
$$
P_{ii}-{\bf P}^w_{ii}=\left(H_{ii}-
{\bf H}_{ii}\right){\bf \Psi}_{n_i}(t_i)W_i^*.
$$
Combining the last equality with \eqref{4.20} gives
\begin{equation}
H_{ii}-{\bf H}_{ii}=\left[\begin{array}{ccc} 0 &
\ldots & 0\\ \vdots & \ddots & \vdots
\\ 0 & \ldots & \delta_i\end{array}\right]\left({\bf
\Psi}_{n_i}(t_i)W_i^*\right)^{-1}.
\label{4.20}
\end{equation}
Since  $|c_{i,0}|=1$, it is seen from definitions \eqref{1.4} and
\eqref{3.3} that the matrix ${\bf \Psi}_{n_i}(t_i)W_i^*$ is upper
triangular and invertible and that its lower diagonal entry equals
\begin{equation}
g_i:=(-1)^{n_i}t_i^{2n_i+1}{c}_{i,0}^*.
\label{4.21}
\end{equation}
Therefore, the inverse matrix $\left({\bf
\Psi}_{n_i}(t_i)W_i^*\right)^{-1}$ is
upper triangular with the lower diagonal entry equal $g_i^{-1}$ so that
the matrix on the right hand side in \eqref{4.20} has all the entries
equal to zero except the lower diagonal entry which is equal to
$\delta_ig_i^{-1}$. Taking into account the definitions \eqref{2.19} and
\eqref{3.3a} we write \eqref{4.20} more explicitly as
$$
\left[c_{i,j+k+1}-w_{j+k+1}(t_i)\right]_{j,k=0}^{n_i}
=\left[\begin{array}{ccc} 0 & \ldots & 0 \\
\vdots & \ddots & \vdots \\ 0 & \ldots &\delta_i
g_{i}^{-1}\end{array}\right].
$$
Upon equating the corresponding entries in the latter equality
we arrive at
$$
w_j(t_i)=c_{i,j}\quad (j=1,\ldots,2n_i)
$$
and
$$
c_{i,2n_i+1}-w_{2n_i+1}(t_i)=\delta_i g_{i}^{-1}.
$$
The first line (together with \eqref{4.16}) proves (\ref{8.8}). The second
one can be written
as
$$
\left(c_{i,2n_i+1}-w_{2n_i+1}(t_i)\right)g_{i}=\delta_i\ge 0,
$$
which implies (\ref{8.9}), due to \eqref{4.21}.
In the case when equality \eqref{4.8a} holds, we get from
(\ref{4.20a}) that $\delta_i=0$ for $i=1,\ldots,k$ and, therefore,
that $w$ is a solution of Problem~$\ref{P:8.2}$ (or
equivalently, of Problem \ref{P:8.1}).\end{proof}

\medskip

We recall now briefly the setting of the Abstract Interpolation
Problem {\bf AIP} (in a generality we need) for the Schur class
$\cS(\cE, \, \cE_*)$ of functions analytic on $\DD$ whose values are
contractive operators mapping a Hilbert space $\cE$ into another
Hilbert space $\cE_*$. The data of the problem consists of Hilbert spaces
$\cE$, $\cE_*$ and $X$, a positive semidefinite linear operator $P$ on
$X$, an operator $T$ on $X$ such that the operator $(I-zT)$ has a
bounded inverse at every point $z\in\overline{\DD}$ except for a
finitely many points, and two linear operators $M: \; X\to \cE$ and
$E: \; X\to \cE_*$ satisfying the identity
\begin{equation}
P-T^*PT=E^*E-M^*M.
\label{4.22}
\end{equation}
\begin{defn}
A function $w\in\cS(\cE, \, \cE_*)$ is said to be a solution of the {\bf
AIP} with the data
\begin{equation}
\{P, \, T, \, E, \, M\} \label{4.23}
\end{equation}
subject to above assumptions, if the function
\begin{equation}
({\bf F}^wx)(t):=\left[\begin{array}{cc} \I_{\cE_*}& w(t) \\ w(t)^*
& \I_{\cE}\end{array}\right]\left[\begin{array}{r} E \\ -M
\end{array}\right]\left(\I-tT\right)^{-1}x.
\label{4.24}
\end{equation}
belongs to the space $H^w$ and
$$
\left\|{\bf F}^w x\right\|_{H^w} \le \|P^{\frac{1}{2}}x\|_X
\quad\mbox{for every}\quad x\in X.
$$
\label{D:8.4}
\end{defn}
The main conclusion of this section is that Problem $\ref{P:6-final}$
can be included into the {\bf AIP} upon specifying the data in
$(\ref{4.22})$ in terms of the data \eqref{data} of Problem
\ref{P:6-final}. Let $X=\C^N$ and $\cE=\cE_*=\C$ and let us
identify the matrices $P$, $T$, $E$ and $M$ defined in
(\ref{3.1})--(\ref{3.7}) with operators acting
between the corresponding finite dimensional spaces. For $T$ of
the form (\ref{3.5}), the operator $(\I-tT)^{-1}$ is well
defined on $X$ for all $t\in\TT\setminus\{t_1,\ldots, t_k\}$.
Also we note that when $X=\C^N$,
$$
\|P^{\frac{1}{2}}x\|_X^2=x^*Px.
$$
Now Theorems \ref{4.2} and \ref{T:4.3} lead us to the following result.
\begin{thm}
Let the matrices $P$, $T$, $E$ and $M$ be given by
$(\ref{3.1})$--$(\ref{3.7})$ and let conditions \eqref{4.9} are satisfied.
Then a Schur function $w$ is a solution of Problem $\ref{P:6-final}$ if
and only if it is a solution of the
{\bf AIP} with the data $(\ref{4.23})$.
\label{equiv-7.7-AIP}
\end{thm}
\begin{cor}
Conditions $P\ge 0$ and $|c_{i,0}|=1$ for $i=1,\ldots,k$ are
necessary and sufficient for Problem $\ref{P:6-final}$ to have a
solution. \label{solution-exists}
\end{cor}
{\bf Proof:} Necessity of the conditions was proved in Theorem
\ref{4.2}. Sufficiency follows from Theorem
\ref{equiv-7.7-AIP} and from a general result \cite{kky}
stating that ${\bf AIP}$ always has a solution.

\section{On the Abstract Interpolation Problem (AIP)}
\setcounter{equation}{0}

In this section we recall some results on the {\bf AIP} formulated
in Definition \ref{D:8.4}. Then in the next section we will
specify these results for the setting of Problem \ref{P:6-final},
when $X=\C^N$, $\cE=\cE_*=\C$ and operators $T$, $E$, $M$ and $P\ge
0$ are just matrices defined in terms of the data of Problem
\ref{P:6-final} via formulas (\ref{3.1})--(\ref{3.7}). In this
section they are assumed to be operators satisfying the Stein
identity (\ref{4.22}) for every $x\in X$. This identity means
that the formula
\begin{equation}
{\bf V}: \; \left[\begin{array}{c}P^{\frac{1}{2}}x \\
Mx\end{array}\right]
\rightarrow \left[\begin{array}{c} P^{\frac{1}{2}}Tx\\
Ex\end{array}\right],\quad x\in X \label{a9.5}
\end{equation}
defines a linear map that  can be extended by continuity to an
isometry ${\bf V}$ acting from
\begin{equation}
{\mathcal {D}}_{\bf V}={\rm Clos}\left\{\left[\begin{array}{c}
P^{\frac{1}{2}}x\\
Mx\end{array}\right], \; x\in X\right\}\subseteq [X]\oplus \cE
\label{a9.6}
\end{equation}
onto
\begin{equation}
{\mathcal {R}}_{\bf V}=
{\rm Clos}\left\{\left[\begin{array}{c}P^{\frac{1}{2}}Tx\\
Ex\end{array}\right], \; x\in X\right\}\subseteq [X]\oplus \cE_*,
\label{a9.7}
\end{equation}
where $[X]={\rm Clos}\{P^{^{\frac{1}{2}}}X\}$. One of the main
results concerning the {\bf AIP} is the characterization of the
set of all solutions in terms of minimal unitary extensions of
${\bf V}$: let ${\mathcal H}$ be a Hilbert spaces containing $[X]$
and let
\begin{equation}
{\bf U}: \; {\mathcal H}\oplus \cE\to {\mathcal H}\oplus \cE_*\quad
({\mathcal H}\supset X) \label{formula10}
\end{equation}
be a unitary operator such that ${\bf U}|_{{\mathcal {D}}_{\bf
V}}={\bf V}$ and having no nonzero reducing subspaces in ${\mathcal
H}\ominus [X]$. Then the {\em characteristic function} of ${\bf
U}$ defined as
\begin{equation}
w(z)={\bf P}_{_{\cE_*}}{\bf U}\left( I-z{\bf P}_{\mathcal H}{\bf
U}\right)^{-1}\vert_{\cE}\quad (z\in\DD), \label{10.300}
\end{equation}
is a solution of the {\bf AIP} and all the solutions to the {\bf AIP} can
be obtained in this way.

A parametrization of all the solutions can be obtained as
follows: introduce the defect spaces
\begin{equation}
\Delta:=\left[\begin{array}{c} [X] \\ \cE\end{array}\right]
\ominus{\mathcal {D}}_{\bf V} \quad{\rm and}\quad
\Delta_*:=\left[\begin{array}{l}[X]\\
\cE_*\end{array}\right]\ominus{\mathcal {R}}_{\bf V} \label{a9.8}
\end{equation}
and let $\widetilde{\Delta}$ and $\widetilde{\Delta}_*$ be
isomorphic copies of $\Delta$ and $\Delta_*$, respectively, with
unitary identification maps
$$
i: \; \Delta\rightarrow \widetilde{\Delta}\quad\mbox{and}\quad
i_*: \; \Delta_*\rightarrow \widetilde{\Delta}_*.
$$
Define a unitary operator ${\bU}_0$ from ${\mathcal {D}}_{\bf
V}\oplus\Delta\oplus\widetilde{\Delta}_*$ onto ${\mathcal
{R}}_{\bf V}\oplus\Delta_*\oplus\widetilde{\Delta}$ by the rule
\begin{equation}
{\bU}_0\vert_{{\mathcal {D}}_{\bf V}}={\bf V},\quad
{\bU}_0\vert_{\Delta}=i,\quad
{\bU}_0\vert_{\widetilde{\Delta}_*}=i_*^{-1}. \label{a9.9}
\end{equation}
This operator is called {\em the universal unitary colligation}
associated to the Stein identity (\ref{4.22}). Since ${\mathcal
{D}}_{\bf V}\oplus\Delta=[X]\oplus \cE$ and ${\mathcal {R}}_{\bf
V}\oplus\Delta_*=[X]\oplus \cE_*$, we can decompose ${\bU}_0$
defined by (\ref{a9.9}) as follows
\begin{equation}
{\bf U}_0=\left[\begin{array}{ccc}U_{11} &U_{12}&U_{13} \\
U_{21}&U_{22}& U_{23}\\ U_{31}&U_{32}&0\end{array}\right]:
\quad\left[\begin{array}{c}[X] \\ \cE \\
\widetilde{\Delta}_*\end{array}
\right]\rightarrow \left[\begin{array}{c}[X] \\ \cE_* \\
\widetilde{\Delta}\end{array}\right]. \label{a9.11}
\end{equation}
Note that $U_{33}=0$, since (by definition (\ref{a9.9})), for
every $\widetilde{\delta}_*\in\widetilde{\Delta}_*$, the vector
${\bU}_0{\delta}_*$ belongs to $\Delta_*$, which is a subspace of
$[X]\oplus \cE_*$ and therefore, is orthogonal to
$\widetilde{\Delta}$. The {\em characteristic function} of
${\bU}_0$ is defined as follows:
\begin{equation}
{\bf S}(z)={\bf P}_{\cE_*\oplus\widetilde{\Delta}}{\bf U}_0\left(
I-z{\bf P}_{[X]}{\bf U}_0\right)^{-1}\vert_{\cE\oplus
\widetilde{\Delta}_*}\quad (z\in\DD), \label{a9.10}
\end{equation}
where ${\bf P}_{\cE_*\oplus\widetilde{\Delta}}$ and ${\bf P}_{[X]}$
are the orthogonal projections of the space $[X]\oplus \cE_*\oplus
\widetilde{\Delta}$ onto $\cE_*\oplus\widetilde{\Delta}$ and $[X]$,
respectively. Upon substituting (\ref{a9.11}) into (\ref{a9.10})
we get a representation of the function ${\bf S}$ in terms of the
block entries of ${\bf U}_0$:
\begin{eqnarray}
{\bf S}(z)&=&\left[\begin{array}{cc}s_0(z) & s_2(z) \\
s_1(z) & s(z)\end{array}\right]\label{a9.11-bis}\\
&=&\left[\begin{array}{cc}U_{22} &U_{23}\\
U_{32}&0\end{array}\right]
+z\left[\begin{array}{c}U_{21}\\U_{31}\end{array}\right]
\left(I_n- zU_{11}\right)^{-1}\left[\begin{array}{cc}U_{12}&
U_{13}\end{array}\right].\nonumber
\end{eqnarray}
The next theorem was proved in \cite{kky}.
\begin{thm}
Let ${\bf S}$ be the characteristic function of the universal
unitary colligation partitioned as in $(\ref{a9.11-bis})$. Then
all the solutions $w$ of the {\bf AIP} are parametrized by the
formula
\begin{equation}
w(z)=s_0(z)+s_2(z)\left(1-\cE(z)s(z)\right)^{-1}\cE(z)s_1(z),
\label{a9.12-bis}
\end{equation}
where $\cE$ runs over the Schur class $\cS(\widetilde{\Delta}, \,
\widetilde{\Delta}_*)$.
\label{T:a9.3}
\end{thm}
Since ${\bf S}$ is the characteristic function of a unitary
colligation, it belongs to the Schur class
$\cS(\cE\oplus\widetilde{\Delta}_*, \;
\cE_*\oplus\widetilde{\Delta})$ (see \cite{Nagy-Foias},
\cite{arovgros1}, \cite{arovgros2}) and therefore,
one can introduce the corresponding de Branges--Rovnyak space
$H^{\bf S}$ as it was explained in Section~2. The next result about
realization of a unitary colligation in a function model space
goes back to M. Livsits, B.~Sz.-Nagy, C.~Foias, L.~de~Branges and J.
Rovnyak. In its present formulation it appears in \cite{kky}--\cite{kh}.

\begin{thm}
Let ${\bf U_0}$ be a unitary colligation of the form
$(\ref{a9.11})$ and let ${\bf S}$ be its characteristic function
defined in $(\ref{a9.10})$. Then the transformation $\cF_{\bU_0}$
defined as
\begin{equation}
\left(\cF_{\bU_0}[x]\right)(z)=
\left[\begin{array}{c}\left(\cF^+_{\bU_0}[x]\right)(z)\\
\left(\cF^-_{\bU_0}[x]\right)(z)
\end{array}\right]:=\left[\begin{array}{c}
{\bf P}_{\cE_*\oplus\widetilde{\Delta}}{\bf U}_0\left(
I-z{\bf P}_{[X]}{\bf U}_0\right)^{-1}[x] \\
\bar{z}{\bf P}_{\cE\oplus\widetilde{\Delta}_*}{\bf U}_0^*\left(
I-\bar{z}{\bf P}_{[X]}{\bf U}_0^*\right)^{-1}[x]
\end{array}\right]
\label{a10.1}
\end{equation}
maps $[X]$ onto the de Branges--Rovnyak space $H^{\bf S}$ and is a
partial isometry.
\label{T:a10.1}
\end{thm}
The transformation $\cF_{\bU_0}$ is called {\em the Fourier
representation} of the space $[X]$ associated with the unitary
colligation ${\bU_0}$ (the name was initiated in
\cite{Nagy-Foias}). Note that the last theorem does not assume any
special structure for ${\bf U}_0$. However, if ${\bU_0}$ is the universal
unitary colligation (\ref{a9.9}) associated to the partially
defined isometry ${\bf V}$ given in (\ref{a9.5}), then
$\cF_{\bU_0}$ can be expressed in terms of $P$, $T$, $E$ and
$M$. The formulation of the following theorem can be found (in a
more general setting) in \cite{khaippars}, \cite{kh}; the proof is
contained in \cite{khphd}. We reproduce it here since the
source is hardly available.
\begin{thm}
Let $\bU_0$ be the universal unitary colligation $(\ref{a9.9})$
associated to the isometry ${\bf V}$ given by $(\ref{a9.5})$ and
let ${\bf S}$ be its characteristic function given by
$(\ref{a9.10})$. Then
\begin{equation}
\left(\cF_{\bU_0}P^{\frac{1}{2}}x\right)(t)=
\left[\begin{array}{cc}I_{\cE_*\oplus\widetilde{\Delta}} & {\bf
S}(t) \\ {\bf S}(t)^* & I_{\cE\oplus \widetilde{\Delta}_*}
\end{array}\right]\left[\begin{array}{c}E(I-tT)^{-1} \\ 0 \\
-M(I-tT)^{-1} \\ 0\end{array}\right]x \label{a10.3}
\end{equation}
for almost every point $t\in\TT$ and for every $x\in X$.
\label{T:a10.2}
\end{thm}
\begin{proof} We will verify (\ref{a10.3}) for ``plus" and ``minus"
components separately, i.e., we will verify the relations
\begin{eqnarray*}
\left(\cF^+_{\bU_0}P^{\frac{1}{2}}x\right)(t)&=&
\left(\left[\begin{array}{c}E\\ 0\end{array}\right]- {\bf
S}(t)\left[\begin{array}{c}M \\ 0
\end{array}\right]\right)(I-tT)^{-1}x,\\
\left(\cF^-_{\bU_0}P^{\frac{1}{2}}x\right)(t)&=& \left({\bf
S}(t)^*\left[\begin{array}{c}E\\ 0\end{array}\right]-
\left[\begin{array}{c}M \\ 0
\end{array}\right]\right)(I-tT)^{-1}x,
\end{eqnarray*}
which are equivalent (upon analytic and conjugate-analytic continuations
inside $\DD$, respectively) to
\begin{eqnarray}
\left(\cF^+_{\bU_0}P^{\frac{1}{2}}x\right)(z)&=&
\left(\left[\begin{array}{c}E\\ 0\end{array}\right]- {\bf
S}(z)\left[\begin{array}{c}M \\ 0
\end{array}\right]\right)(I-zT)^{-1}x,\label{a10.3a}\\
\left(\cF^-_{\bU_0}P^{\frac{1}{2}}x\right)(z)&=& \bar{z}\left({\bf
S}(z)^*\left[\begin{array}{c}E\\ 0\end{array}\right]-
\left[\begin{array}{c}M \\ 0
\end{array}\right]\right)(\bar{z}I-T)^{-1}x. \label{a10.3b}
\end{eqnarray}
To prove (\ref{a10.3a}), we pick an arbitrary vector
$$
v=\left[\begin{array}{c} y \\ e \\
\delta_*\end{array}\right]\in
\left[\begin{array}{c} [X] \\ \cE \\
\widetilde{\Delta}_*\end{array}\right]
$$
and note that by definitions (\ref{a9.10}) and (\ref{a10.1}),
\begin{equation}
{\bf P}_{\cE_*\oplus{\widetilde\Delta}}{\bf U}_0\left(
I-z{\bf P}_{[X]}{\bf U}_0\right)^{-1}\left[\begin{array}{c} y \\ e \\
\delta_*\end{array}\right]=\left(\cF^+_{\bU_0}y\right)(z)+ {\bf
S}(z)\left[\begin{array}{c}e \\ \delta_*\end{array}\right].
\label{a10.4}
\end{equation}
Introduce the vector
$$
v^\prime=\left[\begin{array}{c} y^\prime \\ e^\prime \\
\delta_*^\prime\end{array}\right]:=\left(
I-z{\bf P}_{[X]}{\bf U}_0\right)^{-1}\left[\begin{array}{c} y \\ e \\
\delta_*\end{array}\right]
$$
so that $\left(I-z{\bf P}_{[X]}{\bf U}_0\right)v^\prime=v$.
Comparing the corresponding components in the latter equality
we conclude that $\; e=e^\prime$, $\; \delta_*=\delta_*^\prime$ and
\begin{equation}
y=y^\prime-z{\bf P}_{[X]}{\bf U}_0\left[\begin{array}{c} y^\prime \\
e^\prime \\ \delta_*^\prime\end{array}\right]=
y^\prime-z{\bf P}_{[X]}{\bf U}_0\left[\begin{array}{c} y^\prime \\
e \\ \delta_*\end{array}\right], \label{a10.5}
\end{equation}
so that
\begin{equation}
v^\prime=\left[\begin{array}{c} y^\prime \\ e \\
\delta_*\end{array}\right]=\left(
I-z{\bf P}_{[X]}{\bf U}_0\right)^{-1}\left[\begin{array}{c} y \\ e \\
\delta_*\end{array}\right].
\label{a10.6}
\end{equation}
Substituting (\ref{a10.5}) and (\ref{a10.6}) respectively into
the right and the left hand side expressions in (\ref{a10.4}) we
arrive at
\begin{equation}
{\bf P}_{\cE_*\oplus\widetilde{\Delta}}{\bf U}_0
\left[\begin{array}{c} y^\prime \\ e \\
\delta_*\end{array}\right]=\left(\cF^+_{\bU_0}y^\prime\right)(z)
-z\left(\cF^+_{\bU_0}{\bf P}_{[X]}{\bf U}_0
\left[\begin{array}{c} y^\prime \\ e \\
\delta_*\end{array}\right]\right)(z)+ {\bf
S}(z)\left[\begin{array}{c}e \\ \delta_*\end{array}\right].
\label{a10.7}
\end{equation}
Since the vector $v$ is arbitrary and $I-z{\bf P}_{[X]}{\bf U}_0$
is invertible, it follows by (\ref{a10.5}), that $v^\prime$ can be
chosen arbitrarily in (\ref{a10.7}) . Fix a vector $x\in X$ and
take
\begin{equation}
v^\prime=\left[\begin{array}{c} y^\prime \\ e \\
\delta_*\end{array}\right]=\left[\begin{array}{c}P^{\frac{1}{2}}x \\
Mx \\ 0\end{array}\right]. \label{a10.8}
\end{equation}
Then, by  definition 
(\ref{a9.9}) of ${\bf U}_0$ and definition (\ref{a9.5}) of ${\bf
V}$,
\begin{equation}
{\bU_0}\left[\begin{array}{c}P^{\frac{1}{2}}x \\
Mx\\ 0\end{array}\right]=\left[\begin{array}{c}
P^{\frac{1}{2}}Tx\\ Ex\\ 0\end{array}\right] \label{11.100}
\end{equation}
and thus,
$$
{\bf P}_{[X]}{\bf U}_0\left[\begin{array}{c}P^{\frac{1}{2}}x \\
Mx\\ 0\end{array}\right]=P^{\frac{1}{2}}Tx\quad\mbox{and}\quad
{\bf P}_{\cE_*\oplus\widetilde{\Delta}}{\bf U}_0
\left[\begin{array}{c}P^{\frac{1}{2}}x \\
Mx\\ 0\end{array}\right]=\left[\begin{array}{c}Ex\\
0\end{array}\right].
$$
Plugging the two last relations and (\ref{a10.8}) into
(\ref{a10.7}) we get
$$
\left[\begin{array}{c}Ex\\ 0\end{array}\right]=
\left(\cF^+_{\bU_0}P^{\frac{1}{2}}x\right)(z)-
z\left(\cF^+_{\bU_0}P^{\frac{1}{2}}Tx\right)(z)+ {\bf
S}(z)\left[\begin{array}{c}Mx \\ 0 \end{array}\right].
$$
By linearity of $\cF^+_{\bU_0}$, we have
$$
\left[\begin{array}{c}Ex\\
0\end{array}\right]=\left(\cF^+_{\bU_0}P^{\frac{1}{2}}(I-zT)x\right)(z)
+{\bf S}(z)\left[\begin{array}{c}Mx \\ 0 \end{array}\right]
$$
and, upon replacing $x$ by $(I-zT)x$, we rewrite the last relation
as
$$
\left[\begin{array}{c}E\\
0\end{array}\right](I-zT)^{-1}x=\left(\cF^+_{\bU_0}P^{\frac{1}{2}}x\right)(z)
+{\bf S}(z)\left[\begin{array}{c}M \\ 0
\end{array}\right](I-zT)^{-1}x,
$$
which is equivalent to (\ref{a10.3a}). The proof of (\ref{a10.3b})
is quite similar: we start with an arbitrary vector
$$
v=\left[\begin{array}{c} y \\ e_* \\ \delta\end{array}\right]\in
\left[\begin{array}{c} [X] \\ \cE_* \\
\widetilde{\Delta}\end{array}\right]
$$
and note that by definitions (\ref{a9.10}) and (\ref{a10.1}),
\begin{equation}
\bar{z}{\bf P}_{\cE\oplus\widetilde{\Delta}_*}{\bf U}_0^*\left(
I-\bar{z}{\bf P}_{[X]}{\bf
U}_0^*\right)^{-1}\left[\begin{array}{c} y \\ e_*
\\ \delta\end{array}\right]=\left(\cF^-_{\bU_0}y\right)(z)+
\bar{z}{\bf S}(z)^*\left[\begin{array}{c}e_* \\
\delta\end{array}\right]. \label{a10.9}
\end{equation}
Then we introduce the vector
\begin{equation}
v^\prime:=\left[\begin{array}{c} y^\prime \\ e^\prime_* \\
\delta^\prime\end{array}\right]=\left(
I-\bar{z}{\bf P}_{[X]}{\bf U}_0^*\right)^{-1}\left[\begin{array}{c} y \\ e_* \\
\delta\end{array}\right] \label{a10.10}
\end{equation}
and check that
\begin{equation}
e^\prime_*=e_*,\quad\delta^\prime=\delta,\quad
y=y^\prime-\bar{z}{\bf P}_{[X]}{\bf U}_0^*\left[\begin{array}{c}
y^\prime
\\ e_* \\ \delta\end{array}\right],
\label{a10.10-bis}
\end{equation}
which allows us to rewrite (\ref{a10.9}) as
\begin{equation}
\bar{z}{\bf P}_{\cE\oplus\widetilde{\Delta}_*}{\bf U}_0^*
\left[\begin{array}{c} y^\prime \\ e_* \\
\delta\end{array}\right]=\left(\cF^-_{\bU_0}y^\prime\right)(z)
-\bar{z}\left(\cF^-_{\bU_0}{\bf P}_{[X]}{\bf U}_0^*
\left[\begin{array}{c} y^\prime \\ e_* \\
\delta\end{array}\right]\right)(z)+
\bar{z}{\bf S}(z)^*\left[\begin{array}{c}e_* \\
\delta\end{array}\right]. \label{a10.11}
\end{equation}
By the same arguments as above, $v^\prime$ can be chosen
arbitrarily in $[X]\oplus \cE_*\oplus\widetilde{\Delta}$ and we let
\begin{equation}
v^\prime=\left[\begin{array}{c} y^\prime \\ e \\
\delta_*\end{array}\right]=\left[\begin{array}{c}P^{\frac{1}{2}}Tx \\
Ex \\ 0\end{array}\right],\quad x\in X. \label{a10.12}
\end{equation}
Since ${\bU_0}$ is unitary, it follows from (\ref{11.100}) that
$$
{\bU_0}^*\left[\begin{array}{c}
P^{\frac{1}{2}}Tx\\ Ex\\
0\end{array}\right]=\left[\begin{array}{c}P^{\frac{1}{2}}x \\
Mx\\ 0\end{array}\right]
$$
and thus,
$$
{\bf P}_{[X]}{\bf U}_0^*\left[\begin{array}{c}P^{\frac{1}{2}}Tx \\
Ex\\ 0\end{array}\right]=P^{\frac{1}{2}}x\quad\mbox{and}\quad
{\bf P}_{\cE\oplus\widetilde{\Delta}_*}{\bf U}_0^*
\left[\begin{array}{c}P^{\frac{1}{2}}Tx \\
Ex\\ 0\end{array}\right]=\left[\begin{array}{c}Mx\\
0\end{array}\right].
$$
Plugging the two last relations and (\ref{a10.12}) into
(\ref{a10.11}) we get
$$
\bar{z}\left[\begin{array}{c}Mx\\
0\end{array}\right]=\left(\cF^-_{\bU_0}P^{\frac{1}{2}}Tx\right)(z)
-\bar{z}\left(\cF^-_{\bU_0}P^{\frac{1}{2}}x\right)(z)+ \bar{z}{\bf
S}(z)^*\left[\begin{array}{c}Ex \\ 0\end{array}\right].
$$
By linearity of $\cF^-_{\bU_0}$, we have
$$
\bar{z}\left[\begin{array}{c}Mx\\
0\end{array}\right]=\left(\cF^-_{\bU_0}P^{\frac{1}{2}}(T-\bar{z}I)x\right)(z)
+\bar{z}{\bf S}(z)^*\left[\begin{array}{c}Ex \\
0\end{array}\right]
$$
and, upon replacing $x$ by $(\bar{z}I-T)^{-1}x$, we rewrite the
latter relation as
$$
\bar{z}\left[\begin{array}{c}M\\
0\end{array}\right](\bar{z}I-T)^{-1}=-\left(\cF^-_{\bU_0}P^{\frac{1}{2}}x
\right)(z)+\bar{z}{\bf S}(z)^*\left[\begin{array}{c}E \\
0\end{array}\right] (\bar{z}I-T)^{-1},
$$
which is equivalent to (\ref{a10.3b}).\end{proof}

\section{Description of all solutions of Problem \ref{P:6-final}}
\setcounter{equation}{0}

Since Problem~\ref{P:6-final} is equivalent to the ${\bf AIP}$
with a specific choice of the data (\ref{4.23}), Theorem
\ref{T:a9.3} gives, in fact, a parametrization of all solutions of
Problem~\ref{P:6-final}. However, the fact that in the context of
Problem~\ref{P:6-final}, $X=\C^N$ and $\cE=\cE_*=\C$, and that
the matrices $P$, $T$, $E$ and $M$ are of special structure
\eqref{3.1}-\eqref{3.7}, allow us to rewrite the results from the
previous section more transparently. We assume that the
necessary conditions \eqref{4.9} for Problem  \ref{P:6-final} to have a
solution are in force. Then $P$ satisfies the Stein identity
\begin{equation}
P+M^*M=T^*PT+E^*E \label{5.1}
\end{equation}
(by Theorem \ref{T:8.a1}) which in turn, gives raise to the isometry
$$
{\bf V}: \; \left[\begin{array}{c}P^{\frac{1}{2}}x \\
Mx\end{array}\right]
\rightarrow \left[\begin{array}{c} P^{\frac{1}{2}}Tx\\
Ex\end{array}\right],\quad x\in \C^N
$$
that maps
$$
{\mathcal {D}}_{\bf V}= {\rm Ran}
\, \left[\begin{array}{c}P^{\frac{1}{2}}\\
M\end{array}\right]\subseteq \left[\begin{array}{c}[X] \\
\mathcal{E}\end{array}\right] \quad {\rm onto} \quad {\mathcal
{R}}_{\bf V}={\rm Ran}
\, \left[\begin{array}{c}P^{\frac{1}{2}}T\\
E\end{array}\right]\subseteq \left[\begin{array}{c}[X] \\
\mathcal{E}_*\end{array}\right],
$$
where $[X]= {\rm Ran\ } P^{\frac{1}{2}}$. In the present context,
the defect spaces (\ref{a9.8})
$$
\Delta=\left[\begin{array}{c}[X] \\
\mathcal{E}\end{array}\right] \ominus{\mathcal {D}}_{\bf V}
\quad{\rm and}\quad \Delta_*= \left[\begin{array}{c}[X] \\
\mathcal{E}_*\end{array}\right] \ominus{\mathcal {R}}_{\bf V}
$$
admit a simple characterization.
\begin{lem}
If $P$ is nonsingular, then
\begin{equation}
\Delta={\rm Span}\left[\begin{array}{c}-P^{^{-\frac{1}{2}}}M^*\\
1\end{array}\right]\quad{\rm and}\quad \Delta_*={\rm Span}
\left[\begin{array}{c}-P^{^{-\frac{1}{2}}}(T^{^{-1}})^*E^*\\
1\end{array}\right].
\label{10.60}
\end{equation}
If $P$ is singular then $\Delta=\{0\}$ and $\Delta_*=\{0\}.$
\label{defect-spaces}
\end{lem}
\begin{proof}
A vector $\left[\begin{array}{c} [x] \\
e\end{array}\right]\in \left[\begin{array}{c}[X] \\
\mathcal{E}\end{array}\right]$ belongs to $\Delta$ if and only if
$$<[x], \, P^{\frac{1}{2}}y> + <e, \, My>=0$$ for every $y\in X,$
which is equivalent to
\begin{equation}
P^{\frac{1}{2}}[x] + M^*e=0. \label{10.200}
\end{equation}
Equation (\ref{10.200}) has a nonzero solution
$\left[\begin{array}{c} [x]
\\ e\end{array}\right]$ if and only if the vector-column $M^*$
belongs to $[X]$. If $P$ is nonsingular, then $[X]=X$, therefore
$M^*\in [X]$, and (\ref{10.200}) implies the first relation in
(\ref{10.60}). The second relation is proved quite similarly.

Let now $P$ be singular. Then $M^*\notin [X]$. Indeed assuming
that $M^*\in {\rm Ran} \, P^{\frac{1}{2}}$ we get that $Mx=0$ for
every $x\in{\rm Ker} \, P$, which implies, in view of (\ref{5.1}),
that $Tx\in{\rm Ker} \, P$ and $Ex=0$ for every $x\in{\rm Ker} \,
P$. In particular, ${\rm Ker} \, P$ is $T$-invariant and
therefore, at least one eigenvector $x_0$ of $T$ belongs to ${\rm
Ker} \, P$, and this vector must satisfy  $Ex_0=0$. However, by
definitions (\ref{3.5}, (\ref{3.6}) $Ex_0\ne 0$ for every
eigenvector $x_0$ of $T$. The contradiction means that $M^*\notin
[X]$ and, therefore, equation (\ref{10.200}) has only zero
solution, i.e. $\Delta=\{0\}$ in case when $P$ is singular. The
result concerning $\Delta_*$ is established in much the same way.
\end{proof}
\begin{thm}
If $P$ is singular, then Problem \ref{P:6-final} has a unique solution
\begin{equation}
w_0(z)=E\left( {\widetilde P}-zPT\right)^{-1}M^*,
\label{5.2}
\end{equation}
(which is a finite Blaschke product of degree equal to ${\rm rank}\ P$),
where
\begin{equation}
{\widetilde P}:=P+M^*M=T^*PT+E^*E. \label{5.3}
\end{equation}
The inverse in $(\ref{5.2})$ is well defined as an operator
on $\widetilde X={\rm Ran} \, \widetilde P$.
\label{T:H1}
\end{thm}
\begin{proof} By Lemma \ref{defect-spaces}, if $P$ is singular  then
${\mathcal D}_{\bf V}=[X]\oplus \cE$ and ${\mathcal R}_{\bf
V}=[X]\oplus \cE_*$ where $[X]={\rm Ran} \, P^{\frac{1}{2}}$.
Therefore, the isometry ${\bf V}$ defined by (\ref{a9.5}),
is already a unitary operator from $[X]\oplus \cE$ onto $[X]\oplus
\cE_*$. Therefore, the solution is unique and is given by the formula
(\ref{10.300}) with ${\bf V}$ and $[X]$ in place of  ${\bf U}$ and
${\mathcal H}$, respectively:
\begin{equation}
w_0(z)={\bf P}_{_{\cE_*}}{\bf V}\left( I-z{\bf P}_{[X]}{\bf
V}\right)^{-1}\vert_{\cE}\quad (z\in\DD). \label{10.400}
\end{equation}
Since $\dim[X]<\infty$, it follows that $w_0$ is a finite Blaschke
product of degree equal to $\dim[X]={\rm rank}\
P$ (see, e.g., \cite{Nik}). It remains to derive the realization formula
\eqref{5.2} from (\ref{10.400}).

Note that by definition of ${\widetilde X}$, it is ${\widetilde
P}$-invariant. Since ${\widetilde P}$ is Hermitian, it is invertible on
its range ${\widetilde X}.$ In what follows, the symbol ${\widetilde
P}^{-1}$ will be understood as an operator on ${\widetilde X}$.
We define the mappings
$$
A=\left[\begin{array}{c}P^{\frac{1}{2}}\\
M\end{array}\right]: {\widetilde X}\to [X]\oplus \cE\quad {\rm
and\quad }B=\left[\begin{array}{c}P^{\frac{1}{2}}T\\
E\end{array}\right]: {\widetilde X}\to [X]\oplus \cE_*.
$$
 Since
\begin{equation}
A^*A=B^*B={\widetilde P} \label{5.4}
\end{equation}
and since ${\widetilde P}$ is invertible on ${\widetilde X},$ both
$A$ and $B$ are nonsingular on ${\widetilde X}$.
Since $P$ is singular, it follows (by the proof of Lemma
\ref{defect-spaces}) that $M^*\notin[X]$
and thus $\dim {\widetilde X}=\dim [X]+1$. Therefore, $A$ is a
bijection from ${\widetilde X}$ onto $[X]\oplus \cE$ and $B$ is a
bijection from ${\widetilde X}$ onto $[X]\oplus \cE_*$. Using
(\ref{5.4}), one can also write the formulas for the inverses
$$
A^{-1}={\widetilde P}^{-1}A^*:[X]\oplus \cE\to{\widetilde X} ,\quad
B^{-1}={\widetilde P}^{-1}B^*:[X]\oplus \cE_*\to{\widetilde X}.
$$
By definition (\ref{a9.5}), ${\bf V}A=B,$ which can be rephrased
as ${\bf V}=BA^{-1}=B{\widetilde P}^{-1}A^*$. Plugging this in
(\ref{10.400}) we get \eqref{5.2}:
\begin{eqnarray}
w_0(z) &=& {\bf P}_{_{\cE_*}}B{\widetilde P}^{-1}A^* \left( I-z{\bf
P}_{[X]}B{\widetilde P}^{-1}A^*\right)^{-1}\vert_{\cE}\nonumber\\
&=& {\bf P}_{_{\cE_*}}B{\widetilde P}^{-1} \left( I-zA^*{\bf
P}_{[X]}B{\widetilde P}^{-1}\right)^{-1}A^*\vert_{\cE}\nonumber\\
&=& {\bf P}_{_{\cE_*}}B \left( {\widetilde P}-zA^*{\bf
P}_{[X]}B\right)^{-1}A^*\vert_{\cE}\nonumber\\
&=& {\bf P}_{_{\cE_*}}B \left( {\widetilde
P}-zPT\right)^{-1}A^*\vert_{\cE}\nonumber\\
&=& E\left( {\widetilde P}-zPT\right)^{-1}M^*.\nonumber
\end{eqnarray}
All the inverses in the latter chain of equalities (except the first
one) are understood as operators on ${\widetilde X}.$ They exist,
since the first inverse in this chain does, which, in turn, is in
effect since ${\bf V}$ is unitary.\end{proof}

\begin{thm}
If $P$ is nonsingular, then the set of all solutions of
Problem~\ref{P:6-final} is parametrized by the formula
\begin{equation}
w(z)=s_0(z)+s_2(z)\left(1-\cE(z)s(z)\right)^{-1}\cE(z)s_1(z),
\label{5.5}
\end{equation}
where the free parameter $\cE$ runs over the Schur class $\cS$,
\begin{eqnarray}
s_0(z)&=& E(\widetilde{P}-z P T)^{-1}M^*,\label{5.6}\\
s_1(z)&=&\alpha^{-1}\left(1- z M T(\widetilde{P}-z P
T)^{-1}M^*\right),
\label{5.7}\\
s_2(z)&=&\beta^{-1}\left(1- z E(\widetilde{P}-z P
T)^{-1}(T^{-1})^*E^*\right),
\label{5.8}\\
s(z)&=&z\alpha^{-1}\beta^{-1} M P^{-1}\widetilde P(\widetilde{P}-z
P T)^{-1}(T^{-1})^*E^*,\label{5.9}
\end{eqnarray}
the matrix $\widetilde{P}$ is given in \eqref{5.3}
and $\alpha$ and $\beta$ are positive numbers given by
\begin{equation}
\alpha=\sqrt{1+MP^{-1}M^*}\quad\mbox{and}\quad
\beta=\sqrt{1+ET^{-1}P^{-1}(T^{-1})^*E^*}. \label{5.10}
\end{equation}
The matrix $(\widetilde{P}-zPT)$ is invertible for every $z\in\DD$
in this case.
\label{T:H2}
\end{thm}
\begin{proof} By Theorem \ref{T:a9.3}, all the solutions of Problem
\ref{P:6-final} are parametrized by the formula (\ref{5.5})
where the coefficients $s_0$, $s_1$, $s_2$ and $s$ are the entries
of the characteristic function ${\bf S}$ of the universal unitary
colligation ${\bf U_0}$. By Lemma \ref{defect-spaces}, we have
$\dim \, \Delta=\dim\,\Delta_*=1$. Since, by the very construction
of the universal colligation, $\widetilde{\Delta}$ and
$\widetilde{\Delta}_*$ are isomorphic copies of $\Delta$ and
$\Delta_*$, respectively, we have also $\dim \, \widetilde{\Delta}
=\dim\,\widetilde{\Delta}_*=1$, and we will identify each of these
two spaces with $\C$. However, we will keep the notations
$\widetilde{\Delta}$ and $\widetilde{\Delta}_*$ for the spaces so
that not to mix them up. Thus, in the present context, the
characteristic function ${\bf S}$ of ${\bf U}_0$ is a $2\times 2$
matrix valued function and it remains to establish explicit
formulas (\ref{5.6})--(\ref{5.9}) for its entries
which are scalar valued functions. First we will write relations
(\ref{a9.9}) defining the operator ${\bf U}_0: \; X\oplus \cE\oplus
{\widetilde \Delta}_*\to X\oplus \cE_*\oplus {\widetilde \Delta}$
more explicitly. The first relation in (\ref{a9.9}) can be
written, by the definition (\ref{a9.5}) of ${\bf V}$,  as
\begin{equation}
{\bf U}_0\left[\begin{array}{c} P^{\frac{1}{2}}\\ M \\ 0
\end{array}\right]=\left[\begin{array}{c}P^{\frac{1}{2}}T \\ E \\
0\end{array}\right]. \label{5.11}
\end{equation}
By Lemma \ref{defect-spaces}, the spaces $\Delta$ and $\Delta_*$
are spanned by the vectors
$$
\delta=\left[\begin{array}{c}-P^{^{-\frac{1}{2}}}M^*\\
1 \\ 0\end{array}\right]\in\left[\begin{array}{c}X\\
\cE \\ {\widetilde \Delta}_*\end{array}\right]\quad\mbox{and}\quad
\delta_*=\left[\begin{array}{c}-P^{^{-\frac{1}{2}}}(T^{^{-1}})^*E^*\\
1 \\  0\end{array}\right]\in\left[\begin{array}{c}X\\
\cE_* \\ {\widetilde \Delta}\end{array}\right],
$$
respectively. Note that
\begin{equation}
\|\delta\|^2=1+MP^{-1}M^*\quad\mbox{and}\quad \|\delta_*\|^2
=1+ET^{-1}P^{-1}(T^{-1})^*E^*.
\label{5.12}
\end{equation}
By the second relation in (\ref{a9.9}), the vector ${\bf U}_0\delta$
belongs to $\widetilde{\Delta}$ and therefore, it is of the form
\begin{equation}
{\bf U}_0\delta=\left[\begin{array}{c} 0 \\ 0 \\
\alpha\end{array}\right] \label{5.13}
\end{equation}
where $|\alpha|=\|\delta\|$, due to unitarity of ${\bf U}_0$. The
latter equality and the first equality in (\ref{5.12}) imply that
$\alpha\ne 0$. In fact, we can choose the identification map $i: \;
\Delta\to\widetilde{\Delta}$ so that $\alpha$ will be as in
(\ref{5.10}). Equality (\ref{5.13}) is an
explicit form of the second relation in (\ref{a9.9}). Similarly,
the second identification map $i_*: \;
\Delta_*\to\widetilde{\Delta}_*$ can be chosen so that
\begin{equation}
{\bf U}_0\left[\begin{array}{c} 0 \\ 0 \\ \beta\end{array}\right]
=\delta_*,
\label{5.14}
\end{equation}
where $\beta$ is defined as in (\ref{5.10}). Summarizing equalities
(\ref{5.11}), (\ref{5.13}) and (\ref{5.14}) we conclude
that ${\bf U}_0$ satisfies (and is uniquely determined by) the equation
\begin{equation}
{\bf U}_0A=B, \label{5.15}
\end{equation}
where
\begin{equation}
A=\left[\begin{array}{ccc}
P^{\frac{1}{2}} & -P^{-\frac{1}{2}}M^* & 0\\
M & 1 & 0 \\ 0 & 0 & \beta\end{array}\right]\quad\mbox{and}\quad
B=\left[\begin{array}{ccc}P^{\frac{1}{2}}T & 0 &
-P^{-\frac{1}{2}}(T^{-1})^*E^* \\ E & 0 & 1 \\
0 & \alpha & 0\end{array}\right] \label{5.16}
\end{equation}
are operators from $X\oplus E\oplus \cE_*$ to $X\oplus E\oplus
{\widetilde \Delta}_*$ and to $X\oplus E_*\oplus {\widetilde
\Delta}$, respectively. Since ${\bf U}_0$ is unitary, it follows
that $A^*A=B^*B$. We denote this matrix by $\widehat{P}$ and a
straightforward calculation shows that
\begin{equation}
\widehat{P}:=A^*A=B^*B=\left[\begin{array}{ccc} {\widetilde P} & 0
& 0\\ 0 & |\alpha|^2 & 0 \\ 0 & 0 &
|\beta|^2\end{array}\right]:\quad
\left[\begin{array}{c}X \\ \cE \\
\cE_*\end{array}
\right]\rightarrow \left[\begin{array}{c}X \\ \cE \\
\cE_*\end{array}\right]. \label{5.17}
\end{equation}
where $\widetilde{P}$ is given in (\ref{5.3}). Since $P$ is
nonsingular so is $\widetilde P$ and since $\alpha\ne 0$ and
$\beta\ne 0$, $\widehat P$ is nonsingular as well. Therefore, $A$
and $B$ are nonsingular. Now we proceed as in the proof of Theorem
\ref{T:H1}: it follows from (\ref{5.15}) and (\ref{5.17})
that ${\bf U}_0=BA^{-1}=B{\widehat P}^{-1}A^*$ which being
substituted into (\ref{a9.10}) leads us (recall that since $P$ is
nonsingular, $[X]=X=C^N$) to
\begin{eqnarray}
{\bf S}(z)&=& {\bf P}_{\cE_*\oplus\widetilde{\Delta}}B{\widehat
P}^{-1}A^* \left( I-z{\bf P}_{X}B{\widehat
P}^{-1}A^*\right)^{-1}\vert_{\cE\oplus
\widetilde{\Delta}_*}\nonumber\\
&=& \left[\begin{array}{cc} 0 & I_2\end{array}\right] B{\widehat
P}^{-1} \left( I-zA^*{\bf P}_{X}B{\widehat P}^{-1}\right)^{-1}A^*
\left[\begin{array}{c}0 \\ I_2\end{array}\right]\label{5.18}\\
&=& \left[\begin{array}{ccc} E & 0 & 1\\
0 & \alpha & 0\end{array}\right] \left({\widehat P}-zA^*
\left[\begin{array}{cc}I_N & 0 \\ 0 &
0\end{array}\right]B\right)^{-1} \left[\begin{array}{cc}M^* & 0 \\
1 & 0 \\ 0 & \beta\end{array}\right],\nonumber
\end{eqnarray}
where $I_2$ and $I_N$ are $2\times 2$ and $N\times N$ unit
matrices, respectively. The first inverse in this chain of
equalities exists for every $z\in\DD$ since ${\bf U}_0$ is
unitary, all the others exist since the first one does. By
(\ref{5.16}) and (\ref{5.17}),
$$
\widehat P-zA^*\left[\begin{array}{cc}I_N & 0 \\
0 & 0\end{array}\right]B=\left[\begin{array}{ccc}
{\widetilde P}-zPT & 0 & z(T^{-1})^*E^*\\
zMT & |\alpha|^2 & -zMP^{-1}(T^{-1})^*E^* \\ 0 & 0 &
|\beta|^2\end{array}\right].
$$
Upon inverting the latter triangular matrix and plugging it into
(\ref{5.18}), we eventually get
\begin{eqnarray*}
{\bf S}(z)&=&\left[\begin{array}{cc}s_0(z) & s_2(z) \\
s_1(z) & s(z)\end{array}\right]\\
& &\\
&=&\left[\begin{array}{ccc} ER(z)M^* & &
\beta^{-1}\left(1-zER(z)(T^{-1})^*E^*\right)\\
& &\\ \alpha^{-1}\left(1-zMTR(z)M^*\right) & &
z\alpha^{-1}\beta^{-1} MP^{-1}\widetilde{P}R(z)(T^{-1})^*E^*
\end{array}\right],
\end{eqnarray*}
where $R(z)=\left(\widetilde{P}-zPT\right)^{-1}$, which is
equivalent to (\ref{5.6})--(\ref{5.9}).\end{proof}

In conclusion we will establish some important properties of the
coefficient matrix ${\bf S}$ constructed in Theorem \ref{T:H2}.
\begin{thm}
Let  ${\bf S}=\left[\begin{array}{cc}s_0 & s_2 \\
s_1 & s\end{array}\right]$ be the characteristic function
of the universal unitary colligation $\bU_0$ defined in
\eqref{5.16}, \eqref{5.17}. Then
\begin{enumerate}
\item The function $s_0$ is a solution of Problem \ref{P:8.1}.
\item The function ${\bf S}(z)$ is a rational inner matrix-function
of degree at most $N$.
\item
The functions $s_1$ and $s_2$ have zeroes of multiplicity $n_i+1$
at each interpolating point $t_i$ and do not have other zeroes.
\end{enumerate}
\label{T:a10.3}
\end{thm}
{\bf Proof:} By Theorem \ref{T:H2}, $s_0$ is a solution of Problem
\ref{P:6-final} (corresponding to the parameter $\cE\equiv 0$ in the
parametrization formula (\ref{5.5})). Therefore, by
Theorem~\ref{equiv-7.7-AIP}, ${\bf F}^{s_0}x$ belongs to the space
$H^{s_0}$ for every $x\in\C^N=X$, where $H^{s_0}$ is the de
Branges--Rovnyak space associated to the Schur function $s_0$ and
where
\begin{equation}
{\bf F}^{s_0}(t):=\left[\begin{array}{cc} 1& s_0(t) \\
s_0(t)^* & 1\end{array}\right]\left[\begin{array}{r} E \\ -M
\end{array}\right]\left(\I-tT\right)^{-1}.
\label{5.19}
\end{equation}
Again, by Theorem~\ref{equiv-7.7-AIP}, to show that $s_0$ is a
solution of Problem \ref{P:8.1}, it remains to check that $
\left\|{\bf F}^{s_0}x\right\|_{H^{s_0}}=x^* Px $. Letting
for short
\begin{equation}
R_T(t):=\left(\I-tT\right)^{-1}, \label{5.20}
\end{equation}
we note that by (\ref{5.19}) and by definition of the norm in the
de-Branges-Rovnyak space,
\begin{equation}
\left\|{\bf F}^{s_0}x\right\|^2_{H^{s_0}}=
\left\langle \left[\begin{array}{cc} 1& s_0 \\
s_0^* & 1\end{array}\right]\left[\begin{array}{r} E \\ -M
\end{array}\right]R_Tx, \; \left[\begin{array}{r} E \\ -M
\end{array}\right]R_Tx\right\rangle_{L^2(\C^2)}
\label{5.21}
\end{equation}
By Theorems \ref{T:a10.1} and \ref{T:a10.2}, the function
\begin{eqnarray}
\left(\cF_{\bU_0}P^{\frac{1}{2}}x\right)(t)&=&
\left[\begin{array}{cc}\I_2&
{\bf S}(t) \\ {\bf S}(t)^* & \I_2
\end{array}\right]\left[\begin{array}{r}E \\ 0 \\
-M \\ 0\end{array}\right]R_T(t)x\label{5.22}\\
&=&\left[\begin{array}{c}E-s_0(t)M \\ -s_1(t)M\\
s_0(t)^*E-M \\ s_2(t)^*E \end{array}\right]\left(\I-tT\right)^{-1}x
\quad\mbox{belongs to $H^{\bf S}$}\nonumber
\end{eqnarray}
for every vector $x\in \C^N$. Note that
$E\left(\I-tT\right)^{-1}x\not\equiv 0$, unless $x=0$.
Indeed, letting
\begin{equation}
x=\left[\begin{array}{c}x_1 \\ \vdots \\ x_k\end{array}\right],
\quad\mbox{where}\quad x_i=\left[\begin{array}{c}x_{i,1} \\ \vdots \\
x_{i,n_i}\end{array}\right]\quad(i=1,\ldots,k), \label{5.23}
\end{equation}
we get, on account of definitions (\ref{3.5}) and (\ref{3.6})
of $T$ and $E$, that
\begin{equation}
E\left(\I-tT\right)^{-1}x=\sum_{i=1}^k\sum_{j=0}^{n_i}\frac{t^j}
{(1-t\bar{t}_i)^{j+1}}x_{i,j}\not\equiv 0 \label{5.24}
\end{equation}
for every $x\ne 0$, since the functions
$$
\frac{t^j}{(1-t\bar{t}_i)^{j+1}} \quad (i=1,\ldots,k; \;
j=0,\ldots,n_i)
$$
are linearly independent (recall that all the points
$t_1,\ldots,t_k$ are distinct).

Note also that $s_2\not\equiv 0$ (since $s_2(0)=\beta\neq 0$, by
\eqref{5.8} and \eqref{5.10}) and therefore,
$$
s_2(t)E\left(\I-tT\right)^{-1}x\not\equiv 0
$$
for every $x\in X, x\ne 0$. It is seen from (\ref{5.22}) that the latter
function
is the bottom component of $\cF_{\bU_0}P^{\frac{1}{2}}x$, which
leads us to the conclusion that
$$
\cF_{\bU_0}P^{\frac{1}{2}}x\not\equiv 0\quad\mbox{for every}\quad x\neq 0.
$$
The latter means that the linear map $\cF_{\bU_0}: \; [X]\to
H^{\bf S}$ is a bijection. Since $\cF_{\bU_0}$ is a partial
isometry (by Theorem \ref{T:a10.1}), it now follows that this map
is unitary, i.e., that
\begin{equation}
\left\|\cF_{\bU_0}P^{\frac{1}{2}}x\right\|_{H^{\bf S}}^2=
\left\|P^{\frac{1}{2}}x\right\|_{X}^2=x^*Px. \label{5.25}
\end{equation}
Furthermore,
$$
\left\|\cF_{\bU_0}P^{\frac{1}{2}}x\right\|^2_{H^{\bf S}}=
\left\langle
\left[\begin{array}{cc}\I_2
& {\bf S} \\ {\bf S}^* & \I_2
\end{array}\right]\left[\begin{array}{c}E \\ 0 \\
-M \\ 0\end{array}\right]R_Tx, \; \,
\left[\begin{array}{c}E \\ 0 \\
-M \\ 0\end{array}\right]R_Tx \right\rangle_{L^2(\C^4)},
$$
by (\ref{5.22}) and virtue of formula \eqref{2.1} for the norm in $H^{\bf
S}$. Upon taking advantage of the zero entries in the last formula and the
partition of the matrix ${\bf S}$, we get
\begin{eqnarray}
\left\|\cF_{\bU_0}P^{\frac{1}{2}}x\right\|^2_{H^{\bf S}}
&=&\left\langle\left[\begin{array}{cccc}1
& 0 & s_0 & s_2 \\ s_0^* & s_1^* & 1  & 0 \end{array}\right]
\left[\begin{array}{c}E \\ 0 \\
-M \\ 0\end{array}\right]R_Tx, \; \, \left[\begin{array}{c}E \\
-M\end{array}\right]R_Tx
\right\rangle_{L^2(\C^2)}\nonumber\\
&=&\left\langle \left[\begin{array}{cc} 1
& s_0 \\ s_0^* & 1
\end{array}\right]\left[\begin{array}{r} E \\
-M \end{array}\right]R_Tx, \; \left[\begin{array}{r} E \\ -M
\end{array}\right]R_Tx\right\rangle_{L^2(\C^2)}
\label{5.26}
\end{eqnarray}
Comparing  (\ref{5.21}) and (\ref{5.26}) and taking into account
(\ref{5.25}) we come to
$$
\left\|{\bf F}^{s_0}x\right\|_{H^{s_0}}=
\left\|\cF_{\bU_0}P^{\frac{1}{2}}x\right\|_{H^{\bf S}}=x^*Px,
$$
which proves the first assertion of the theorem.
The second assertion 
follows since $\dim X=N={\displaystyle \sum_{i=0}^k(n_i+1)}<\infty
$ (see, e.g., \cite{Nik}).

To prove the last assertion, we use (\ref{5.22}) for $x$ in the form
(\ref{5.23}) with the only nonzero entry $x_{i,n_i}=1$. For this choice of
$x$ we have by definitions \eqref{3.5}--\eqref{3.7} of $T$, $E$ and $N$,
$$
E (I-tT)^{-1}x=\frac{t^{n_i}}{(1-t\bar{t}_i)^{n_i+1}}
\quad\mbox{and}\quad
M(I-tT)^{-1}x=\frac{{\bf c}_i(t)}{(1-t\bar{t}_i)^{n_i+1}}
$$
where
\begin{equation}
{\bf c}_i(t)={\displaystyle\sum_{\ell=0}^{n_i}t^{n_i-\ell}
(1-t\bar{t}_i)^{\ell}c_{i,\ell}^*}.
\label{5.27}
\end{equation}
Now we conclude from \eqref{5.22} that
\begin{equation}
\frac{s_1(t){\bf c}_i(t)}{(1-t\bar{t}_i)^{n_i+1}}\in H_2^+
\quad\mbox{and}\quad \frac{t^{n_i}s_2(t)^*}{(1-t\bar{t}_i)^{n_i+1}}=
\bar{t}\frac{s_2(t)^*}{(\bar{t}-\bar{t}_i)^{n_i+1}}\in H_2^-.
\label{5.28}
\end{equation}
By (\ref{5.27}), ${\bf c}_i(t_i)=t_i^{n_i}c_{i,0}^*\neq 0$ and thus,
the first condition
in (\ref{5.28}) implies that $s_1$ has the zero of multiplicity
at least $n_i+1$ at $t_i$. The second condition in (\ref{5.28})
is equivalent to
$$
\frac{s_2(t)}{(t-t_i)^{n_i+1}}\in H_2^+
$$
which implies that $s_2$ has zero of multiplicity at least $n_i+1$
at $t_i$. On the other hand, since $s_1$ and $s_2$ are rational
functions of degree at most $N={\displaystyle\sum_{i=0}^k(n_i+1)}$
(the second assertion of this theorem) and since they do not
vanish identically (by the proof of the first assertion of this
theorem), they can not have more than $N$ zeroes. Therefore, they
have zeroes of multiplicities $n_i+1$ at $t_i$ for $i=1,\ldots,k$
and they do not have other zeroes. \qed

\bigskip

Some consequences of Theorem \ref{T:a10.3} needed in the next
section are proved in the following lemma.
\begin{lem} Let ${\bf S}=\left[\begin{array}{cc}s_0 & s_2 \\
s_1 & s\end{array}\right]$ be as in Theorem $\ref{T:a10.3}$. Then
$|s(t_i)|=1$,
\begin{equation}
\frac{s_1^{(n_i+1)}(t_i)}{(n_i+1)!}=\lim_{z\to
t_i}\frac{s_1(z)}{(z-{t}_i)^{n_i+1}}\neq 0,\quad
\frac{s_2^{(n_i+1)}(t_i)}{(n_i+1)!}=\lim_{z\to
t_i}\frac{s_2(z)}{(z-{t}_i)^{n_i+1}}\neq 0,  \label{5.29}
\end{equation}
and
\begin{equation}
{s_2^{(n_i+1)}(t_i})^*=(-1)^{n_i} t_i^{2n_i+2}
{s(t_i})^*s_1^{(n_i+1)}(t_i) {c_{i,0}}^*\ .
\label{5.30}\end{equation}
\label{L:6.5}
\end{lem}
{\bf Proof:} By the third assertion of Theorem \ref{T:a10.3}, the
rational functions $s_1$ and $s_2$ have zeros of multiplicity
$n_i+1$ at $t_i$. This implies  (\ref{5.29}). By the second assertion of
Theorem \ref{T:a10.3}, the matrix-function
${\bf S}$ is inner and rational. In particular, it is unitary at
$t_i\in\T$ and therefore, $|s_2(t_i)|^2+|s(t_i)|^2=1$ which implies
$|s(t_i)|=1$, since $s_2(t_i)=0$. Furthermore, by the reflection
principle, ${\bf S}(1/\bar{z})^*{\bf S}(z)\equiv {\rm I}_2$,or in more
detail,
$$
\left[\begin{array}{cc}s_0(1/\bar{z})^* & s_1(1/\bar{z})^*
\\ s_2(1/\bar{z})^* & s(1/\bar{z})^*\end{array}\right]
\left[\begin{array}{cc}s_0(z) & s_2(z) \\
s_1(z) & s(z)\end{array}\right]=
\left[\begin{array}{cc}1 & 0 \\ 0 & 1\end{array}\right].
$$
In particular,
\begin{equation}
s_2(1/\bar{z})^*s_0(z)+s(1/\bar{z})^*s_1(z)=0 \label{5.31}
\end{equation}
To verify (\ref{5.30}), we note first that by the first relation in
(\ref{5.29}),
\begin{equation}
\lim_{z\to t_i}\frac{s(1/\bar{z})^*s_1(z)}{(z-{t}_i)^{n_i+1}}
=\frac{s(t_i)^*s_1^{(n_i+1)}(t_i)}{(n_i+1)!}. \label{opa}
\end{equation}
Since $|t_i|=1$, the second relation in (\ref{5.29}) gives
$$
\frac{s_2^{(n_i+1)}(t_i)}{(n_i+1)!}=\lim_{z\to t_i}
\frac{s_2(\bar{z}^{-1})}{(\bar{z}^{-1}-{t}_i)^{n_i+1}},
$$
which is equivalent, on account of
$\; \bar{z}^{-1}-{t}_i=-\frac{\bar{z}-\bar{t}_i}{\bar{z}\bar{t}_i}$, to
$$
\frac{s_2^{(n_i+1)}(t_i)}{(n_i+1)!}=
\lim_{z\to t_i}\frac{(-\bar{z}\bar{t}_i)^{n_i+1}
s_2(\bar{z}^{-1})}{(\bar{z}-\bar{t}_i)^{n_i+1}}=(-1)^{n_i+1}\bar{t}_i^{2n_i+2}\lim_{z\to
t_i}\frac{
s_2(\bar{z}^{-1})}{(\bar{z}-\bar{t}_i)^{n_i+1}}.
$$
Upon taking adjoints in the latter equality we get
$$
\lim_{z\to t_i}\frac{s_2(\bar{z}^{-1})^*}{({z}-{t}_i)^{n_i+1}}
=(-1)^{n_i+1}\bar{t}_i^{2n_i+2}\frac{s_2^{(n_i+1)}(t_i)^*}{(n_i+1)!}
$$
and, since $s_0(t_i)=c_{i,0}$ (recall that $s_0$ is a solution of Problem
\ref{P:8.1}), we have also
\begin{equation}
\lim_{z\to t_i}\frac{s_2(\bar{z}^{-1})^*s_0(z)}{({z}-{t}_i)^{n_i+1}}
=(-1)^{n_i+1}\bar{t}_i^{2n_i+2}\frac{s_2^{(n_i+1)}(t_i)^*}{(n_i+1)!}c_{i,0}.
\label{5.32}
\end{equation}
Now upon multiplying (\ref{5.31}) by
${\displaystyle\frac{(n_i+1)!}{(z-t_i)^{n_i+1}}}$ and passing to
limits as $z\to t_i$, we come, on account of (\ref{opa}) and
(\ref{5.32}) to the equality
$$
s(t_i)^*s_1^{(n_i+1)}(t_i)+
(-1)^{n_i+1}\bar{t}_i^{2n_i+2}s_2^{(n_i+1)}(t_i)^*c_{i,0}=0,
$$
which is equivalent to (\ref{5.30}), since $|c_{i,0}|=|t_i|=1$.\qed

\section{Boundary interpolation problem with equality}
\setcounter{equation}{0}

In this section we establish a parametrization of all solutions of
Problem \ref{P:8.1}. Recall that all solutions $w$ of Problem
\ref{P:6-final} are parametrized by the linear fractional formula
(\ref{5.5}) with the free  Schur class parameter $\cE$. Thus,
for every function $w$ of the form (\ref{5.5}), we have
$$
\delta_{w,i}:=\gamma_i-d_{w,n_i}(t_i)\ge 0 \quad(i=1,\ldots,k).
$$
Theorem \ref{T:6.1} below will present the explicit formula for the gaps
$\delta_{w,i}$ in terms of the parameter $\cE$ leading to
$w$ via formula (\ref{5.5}) . As a consequence of this formula we will get
a characterization of
all the parameters $\cE$, leading to functions $w$ with zero gaps,
i.e., to solutions of Problem \ref{P:8.1}. We start with some needed
preliminaries. The proof of the first lemma  can be found in
\cite{sarasonsubh} for the case when $n=0$. For the case $n>0$ the proof
was given in \cite{boldym1} using pretty much the
same ideas.
\begin{lem}
Let $w$ be a function analytic in some nontangential neighborhood
of a point $t_0\in\TT$ and let $w_0,\ldots,w_{2n+1}$ be complex
numbers. Then equality
$$
\lim_{z\to t_0}\frac{w(z)-w_0-(z-t_0)w_1-\ldots -
(z-t_0)^{2n}w_{2n}}{(z-t_0)^{2n+1}}=w_{2n+1}
$$
holds if and only if the nontangential limits
${\displaystyle\lim_{z\to t_0}\frac{w^{(j)}(z)}{j!}}$ exist and
equal $w_j$ for $j=0,\ldots,2n+~1$. \label{L:0.5}
\end{lem}
With every triple $(\omega,t_0,b)$ consisting of a Schur function
$\omega\in\cS$, of a point $t_0\in\T$ and a number $b\in\C$, we
associate the quantity
\begin{eqnarray}
D_{\omega,b}(t_0)
&:=&\int_{\T}\frac{1}{|1-t {\overline
t}_0|^2}\left[\begin{array}{cc}1 &
-b\end{array}\right]\left[\begin{array}{cc}1 & \omega(t) \\
\omega(t)^* & 1\end{array}\right]\left[\begin{array}{c}1 \\
-b^*\end{array}\right] \, m(dt)\nonumber\\
&=&\int_{\T}\left(\ \left\vert \frac{1-\omega(t)\overline
b}{1-t\overline
t_0}\right\vert^2+|b|^2\frac{1-|\omega(t)|^2}{|1-t\overline
t_0|^2} \ \right)m(dt)
\label{7.1}\\
&=&\int_{\T}\left(\ \left\vert
\frac{\omega(t)-b}{t-t_0}\right\vert^2+\frac{1-|\omega(t)|^2}{|t-t_0|^2}\
\right)m(dt),\nonumber
\end{eqnarray}
where $m(dt)$ is the normalized Lebesgue measure on $\T$. It follows from
the very definition that
$$
0\le D_{\omega,\, b}(t_0)\le \infty.
$$
The next theorem (which is a variation of the classical Julia-Carath\'
eodory Theorem and can be mostly found in \cite{sarasonsubh})
characterizes the cases when $D_{\omega,b}(t_0)$ is zero, positive or
infinite. It also establishes relation (\ref{7.3}) which holds true
for any triple $\{\omega, \, t_0, \, b\}$ with $|b|\le 1$, regardless
the case.
\begin{thm} Let $\omega\in\cS$, $t_0\in\T$, $b\in\C$ and let $D_{\omega,\,
b}(t_0)$
be defined as in $(\ref{7.1})$. Then
\begin{enumerate}
\item $D_{\omega,\ b}(t_0)<\infty$ if and only if
\begin{equation}
\liminf_{z\to t_0} \frac{1-|\omega(z)|^2}{1-|z|^2}<\infty \quad
{\rm and}\quad \lim_{z\to t_0}\omega(z)=b,
\label{7.2}
\end{equation}
where the second limit is understood as nontangential. In this
case $|b|=1.$
\item $D_{\omega,\ b}(t_0)=\infty$ if and only if
either
$$
\liminf_{z\to t_0} \frac{1-|\omega(z)|^2}{1-|z|^2}=\infty ,
$$
or the function $\omega$ fails to have a nontangential limit $b$
at $t_0.$
\item $D_{\omega,\
b}(t_0)=0$ if and only if $\omega (z)\equiv b$ and $|b|=1.$
\item If $|b|\le 1,$ then the following equality
\begin{equation}
\lim_{z\to t_0} \frac{1-\omega(z)b^*}{1-z{\overline
t}_0}=D_{\omega,\ b}(t_0)
 \label{7.3}
\end{equation}
holds for the nontangential limit in any event (regardless whether
$D_{\omega,\ b}(t_0)$ is finite or infinite).
\end{enumerate}
\label{L:7.1}
\end{thm}
\begin{proof} Let $H^{\omega}$ be the de Branges-Rovnyak space
associated to the Schur class function $\omega$ and let us consider the
function
\begin{equation}
K_{t_0,b}(t)=\left[\begin{array}{cc}1 &
\omega(t)\\ \omega(t) ^* & 1\end{array}\right] \left[\begin{array}{c}
1\\ -b^*\end{array}\right]\frac {1}{1-t{\overline t}_0}
=\left[\begin{array}{c}K_{t_0,b,+}(t) \\ K_{t_0,b,-}(t)\end{array}\right]
\label{7.4}
\end{equation}
where
\begin{equation}
K_{t_0,b,+}(t)=\frac{1-\omega(t)b^*}{1-t\overline
t_0}\quad\mbox{and}
\quad K_{t_0,b,-}(t)=\overline t\ \frac{\overline{\omega(t})-
b^*}{\overline t-\overline t_0}
\label{7.5}
\end{equation}
By the formula \eqref{2.1},
$$
\|K_{t_0,b}\|^2_{H^{\omega}}=\left\langle
\left[\begin{array}{cc}1 &
\omega(t)\\ \omega(t) ^* & 1\end{array}\right] \left[\begin{array}{c}
1\\ -b^*\end{array}\right]\frac {1}{1-t{\overline t}_0},
\; \left[\begin{array}{c}
1\\ -b^*\end{array}\right]\frac {1}{1-t{\overline
t}_0}\right\rangle_{L^2\oplus L^2}
$$
which is equal to the first integral in \eqref{7.1}.
Therefore, $D_{\omega,\ b}(t_0)<\infty$ if and only if $K_{t_0,b}$
belongs to $L^\omega$ and in this case,
\begin{equation}
D_{\omega,\ b}(t_0)=\Vert K_{t_0,\ b}\Vert_{L^\omega}^2.
\label{7.6}
\end{equation}
On the other hand, if $D_{\omega,\, b}(t_0)<\infty$, then it
follows from the second form of $D_{\omega,\ b}(t_0)$ in
(\ref{7.1}) that
$$
\int_{\T}\left\vert \frac{1-\omega(t)b^*}{1-t\overline{t_0}}
\right\vert^2 \ m(dt)<\infty, \quad\mbox{i.e., that}
\quad K_{t_0,b,+}(t)=\frac{1-\omega(t)b^*}{1-t\overline t_0}\in
L_2.
$$
Since $1-t{\overline t}_0$ is an outer function, it follows, by
Smirnov maximum principle \cite{sm}, that  $K_{t_0,b,+}\in H_2^+$.
Similarly, it follows from the third representation of
$D_{\omega,\ b}(t_0)$ in (\ref{7.1}) that $K_{t_0,b,-}\in H_2^-$.
Therefore, $K_{t_0,b}$ belongs to $H^{\omega}$ by Definition \ref{D:2.1}.
Thus, we have shown that
$$
D_{\omega,b}(t_0)<\infty \; \; \Longleftrightarrow \; \;
K_{t_0,b}\in L^\omega \; \; \Longleftrightarrow \; \; K_{t_0, b}\in
H^\omega.
$$
Now the first assertion of the lemma follows from Theorem
\ref{T:3.2} (the case when $n=0$): the function $K_{t_0,b}$ of
the form (\ref{7.4}) belongs to $H^\omega$ if and only if
conditions in (\ref{7.2}) are satisfied. In this case $|b|=1$,
since
$$
1-|b|^2=\lim_{z\to t_0}(1-|\omega(z)|^2)=\lim_{z\to
t_0}\frac{1-|\omega(z)|^2}{1-|z|^2}(1-|z|^2)=0.
$$
The second assertion is simply the formal negation of the first
one. To prove the third assertion, we observe that $D_{\omega,\ b}(t_0)=0$
if and only if
$$
\left[\begin{array}{cc}1 &
-b\end{array}\right]\left[\begin{array}{cc}1 & \omega(t) \\
\omega(t)^* & 1\end{array}\right]\left[\begin{array}{c}1 \\
-b^*\end{array}\right]=0
$$
almost everywhere on $\TT,$ which occurs if and only if
$$
\left[\begin{array}{cc}1 &
-b\end{array}\right]\left[\begin{array}{cc}1 & \omega(t) \\
\omega(t)^* & 1\end{array}\right]=0
$$
almost everywhere on $\TT$. The latter equality collapses to
$\omega(t)-b=1-b\omega(t)^*=0$ which implies the requisite.

The proof of the fourth assertion splits up into three cases.

\medskip
\noindent
{\bf Case 1}: Let $D_{\omega,\, b}(t_0)<\infty$. Then by the first
statement, conditions \eqref{7.2} are satisfied. Then by Theorem
\ref{T:3.1}, the kernels
$$
K_{z}(t)=\left[\begin{array}{cc}1 & \omega(t)\\ \omega(t) ^* &
1\end{array}\right] \left[\begin{array}{c} 1\\
-\omega(z)^*\end{array}\right]\frac {1}{1-t{\overline z}}
$$
converge to $K_{t_0,\, b}$ in norm of $H^\omega$:
\begin{equation}
K_z\stackrel{\ \ H^w}{\longrightarrow}K_{t_0,\, b},
\label{7.7}
\end{equation}
as $z\to t_0$ nontangentially. By the reproducing property \eqref{2.3}
(for $j=0$),
\begin{equation}
\left\langle f, \; K_{z}\right\rangle_{H^\omega}=f_+(z)\quad
\mbox{for every} \quad f=\left[\begin{array}{c}f_+ \\
f_-\end{array}\right]\in H^\omega.
\label{7.8}
\end{equation}
Then, upon making subsequent use of (\ref{7.6}),
(\ref{7.7}), (\ref{7.8}) and of the explicit formula \eqref{7.5} for
$K_{t_0,b,+}$, we get \eqref{7.3}:
\begin{eqnarray}
D_{\omega,\, b}(t_0)=\Vert K_{t_0,\,
b}\Vert_{H^\omega}^2
&=&\lim_{z\to t_0}\langle K_{t_0,\, b}, K_z\rangle_{H^\omega}\nonumber\\
&=&\lim_{z\to t_0}K_{t_0,b,+}(z) =\lim_{z\to
t_0}\frac{1-\omega(z)b^*}{1-z{\overline t}_0}\ .
\label{7.9}
\end{eqnarray}
{\bf Case 2}: Let $D_{\omega,\, b}(t_0)=\infty$ and
${\displaystyle\liminf_{z\to t_0}
\frac{1-|\omega(z)|^2}{1-|z|^2}<\infty}$.

\medskip

The second assumption guarantees (by Theorem \ref{T:3.2}), that
there exists the nontangential limit $\omega(t_0)={\displaystyle\lim_{z\to
t_0}}\omega(z)$ and that the function  $K_{t_0, \omega(t_0)}$ defined via
(\ref{7.4}), belongs to $H^\omega$. Then  by virtue of (\ref{7.9}), we
have
$$
{\displaystyle \lim_{z\to t_0}} \frac{1-\omega(z)
\omega(t_0)^*}{1-z{\overline t}_0}=D_{\omega,
\omega(t_0)}(t_0)=\|K_{t_0,\omega(t_0)}\|^2_{H^\omega}<\infty.
$$
Since $D_{\omega,b}(t_0)=\infty$, it follows that $b\ne
\omega(t_0)$. It remains to note that \eqref{7.3} again holds since
$$
\lim_{z\to t_0} \frac{1-\omega(z)b^*}{1-z{\overline
t}_0}= \lim_{z\to t_0} \frac{1-\omega(z)\omega(t_0)^*
+\omega(z)(b^*-\omega(t_0)^*)}{1-z{\overline t}_0}=\infty.
$$
{\bf Case 3}: Let $D_{\omega,\, b}(t_0)=\infty$ and
${\displaystyle\liminf_{z\to t_0}
\frac{1-|\omega(z)|^2}{1-|z|^2}=\infty}$. Since
\begin{eqnarray*}
2\Re(1-\omega(z)b^*&=&(1-\omega(z)b^*)+(1-b\omega(z)^*)\\
&=& |1-\omega(z)b^*|^2+1-|b|^2|\omega(z)|^2
\ge 1-|b|^2|\omega(z)|^2,
\end{eqnarray*}
it follows that if $|b|\le 1$, then
\begin{equation}
|1-\omega(z)b^*|\ge \Re(1-\omega(z)b^*)\ge
\frac{1}{2}(1-|\omega(z)|^2). \label{7.10}
\end{equation}
Furthermore, for every $z$ in the following nontangential
neighborhood
$$
\Gamma_a(t_0)= \left\{z\in\DD: \; |t_0-z|<a(1-|z|)\right\},\quad
a>1,
$$
of $t_0$, we have
$$
\frac{1-|z|^2}{|1-z\bar{t}_0|}\ge\frac{1-|z|}{|1-z\bar{t}_0|}>\frac{1}{a}
$$
which together with (\ref{7.10}) leads us to
$$
\left|\frac{1-\omega(z)b^*}{1-z{\overline
t}_0}\right| \ge \frac{1}{2} \,
\frac{1-|\omega(z)|^2}{|1-z\bar{t}_0|}= \frac{1}{2} \,
\frac{1-|\omega(z)|^2}{1-|z|^2} \, \cdot \, \frac{1-|z|^2}
{|1-z\bar{t}_0|}>\frac{1}{2a} \, \frac{1-|\omega(z)|^2}{1-|z|^2}.
$$
Therefore,
$$
\lim_{z\to t_0}\frac{1-\omega(z)b^*}{1-z{\overline
t}_0} =\infty=D_{\omega,\ {b}}(t_0),
$$
which completes the proof of the theorem. \end{proof}
\begin{cor}
If a Schur function $\omega$ is analytic in a neighborhood of
$t_0\in\T$ and $| \omega (t_0)|=~1$, then $D_{\omega,\,
\omega(t_0)}(t_0)<\infty$. In particular, $D_{\omega,\,
\omega(t_0)}(t_0)<\infty$ for every rational $\omega\in\cS$  with $|
\omega (t_0)|=1$. \label{C:12.1}
\end{cor}
\begin{proof} If $w$ meets the assumed properties, then the limit
$$
\lim_{z\to t_0} \frac{1-\omega(z){\overline {
\omega(t_0})}}{1-z{\overline t}_0}={\lim_{z\to t_0}}\left(
\frac{\omega(t_0)-\omega(z)}{t_0-z}\right)\frac{\overline
{\omega(t_0)}} {\overline t_0}=\omega^\prime(t_0) \frac{\overline
{\omega(t_0)}} {\overline t_0}
$$
is finite, then, by the fourth assertion in Lemma \ref{L:7.1},
$D_{\omega,\, {\omega(t_0)}}(t_0)<\infty$.\end{proof}

The next theorem presents an explicit formula for the gap
$\gamma_i-d_{w,n_i}(t_i)$ for any solution $w$ of Problem
\ref{P:6-final}. Recall that by Theorem 
\ref{T:H2}, all solutions of Problem \ref{P:6-final} are
parametrized by formula  (\ref{5.5}).
\begin{thm}
Let $w$ be a solution of Problem $\ref{P:6-final}$, i.e. a
function of the form $(\ref{5.5})$
\begin{equation}
w(z)=s_0(z)+s_2(z)\left(1-\cE(z)s(z)\right)^{-1}\cE(z)s_1(z)
\label{7.11}
\end{equation}
with a parameter $\cE\in{\mathcal S}$. Then for $i=1,\ldots,k$,
\begin{equation}
\gamma_i-d_{w,n_i}(t_i)=\frac{1}{((n_i+1)!)^2} \, \cdot \,
\frac{|s_2^{(n_i+1)}(t_i)|^2}
{D_{\cE,{s(t_i})^*}(t_i)+D_{s,s(t_i)}(t_i)} \label{7.12}
\end{equation}
for $i=1,\ldots,k$, where $D_{\cE,{s(t_i})^*}(t_i)$ and
$D_{s,s(t_i)}(t_i)$ are defined according to $(\ref{7.1})$.
\label{T:6.1}
\end{thm}
{\bf Proof:} Since $w$ is a solution of Problem \ref{P:6-final}
and therefore satisfies conditions \eqref{8.7}--\eqref{8.9},
it follows by Lemma \ref{L:0.5} that
$$
w_{2n_i+1}(t_i)=\lim_{z\to
t_i}\frac{w(z)-c_{i,0}-(z-t_i)c_{i,1}-\ldots
-(z-t_i)^{2n}c_{i,2n_i}}{(z-t_i)^{2n_i+1}}
$$
for $i=1,\ldots,k$. Since $s_0$ is a solution of Problem
\ref{P:8.1} (by the first statement in Theorem \ref{T:a10.3}),
we have (again by Lemma \ref{L:0.5})
$$
c_{i,2n_i+1}=\lim_{z\to t_i}\frac{s_0(z)-c_{i,0}-(z-t_i)c_{i,1}-\ldots
-(z-t_i)^{2n}c_{i,2n_i}}{(z-t_i)^{2n_i+1}}.
$$
Now it follows from the two latter equalities that
$$
c_{i,2n_i+1}-w_{2n_i+1}(t_i)=\lim_{z\to
t_i}\frac{s_0(z)-w(z)}{(z-t_i)^{2n_i+1}},
$$
which being substituted into (\ref{8.7w}), leads us to to
$$
\gamma_i-d_{w,n_i}(t_i)=(-1)^{n_i}t_i^{2n_i+1} \lim_{z\to
t_i}\frac{s_0(z)-w(z)}{(z-t_i)^{2n_i+1}} \, {c_{i,0}}^*.
$$
Substituting
(\ref{7.11}) into the latter equality gives
\begin{equation}
\gamma_i-d_{w,n_i}(t_i)=-(-1)^{n_i}t_i^{2n_i+1} \lim_{z\to
t_i}\frac{s_2(z)(1-\cE(z)s(z))^{-1}\cE(z)s_1(z)}
{(z-t_i)^{2n_i+1}} \; {c^*_{i,0}}. \label{6.2}
\end{equation}
Taking into account relations (\ref{5.29}) (i.e., the fact that
$t_i$ is a zero of multiplicity $n_i$ of $s_1$ and $s_2$), we
rephrase (\ref{6.2}) as
\begin{equation}
\gamma_i-d_{w,n_i}(t_i)=\frac{(-1)^{n_i} t_i^{2n_i+2}}
{((n_i+1)!)^2}s_2^{(n_i+1)}(t_i) \, \lim_{z\to
t_i}\frac{(1-z\bar{t}_i)\cE(z)}{1-\cE(z)s(z)}s_1^{(n_i+1)}(t_i){c^*_{i,0}}.
\label{6.4}
\end{equation}
Due to (\ref{5.30}), the latter equality simplifies to
\begin{equation}
\gamma_i-d_{w,n_i}(t_i)=\frac{|s_2^{(n_i+1)}(t_i)|^2}{((n_i+1)!)^2}
\, \lim_{z\to
t_i}\frac{\cE(z)s(t_i)}{\displaystyle\frac{1-\cE(z)s(z)}{1-z\bar{t}_i}}.
\label{6.10}
\end{equation}
Since $|s(t_i)|=1$ (by Lemma \ref{L:6.5}), we have
\begin{equation}
\frac{1-\cE(z)s(z)}{1-z\bar{t}_i}=\frac{1-s(z){s(t_i})^*}{1-z\bar{t}_i}+
s(z)\frac{1-\cE(z)s(t_i)}{1-z\bar{t}_i}{s(t_i})^*. \label{12.2000}
\end{equation}
 By the fourth assertion of Lemma \ref{L:7.1},
\begin{equation}
\lim_{z\to t_i}\frac{1-s(z){s(t_i})^*}
{1-z\bar{t}_i}=D_{s,s(t_i)}(t_i),\quad \lim_{z\to
t_i}\frac{1-\cE(z)s(t_i)}{1-z\bar{t}_i}= D_{\cE,{s(t_i})^*}(t_i).
\label{12.2100}
\end{equation}
Taking advantage of (\ref{12.2100}) we pass to limits in
(\ref{12.2000}) as $z\to t_i$ to get
\begin{equation}
\lim_{z\to t_i}\frac{1-\cE(z)s(z)}{1-z\bar{t}_i}=
D_{s,s(t_i)}(t_i)+D_{\cE,{s(t_i})^*}(t_i). \label{12.2200}
\end{equation}
Since $s$ is rational and $|s(t_i)|=1$, it follows by
Corollary \ref{C:12.1} that $D_{s,s(t_i)}(t_i)$ is
finite. Since, by Theorem \ref{T:a10.3}, ${\bf S}(t)$ is unitary
for $t\in\TT$ and $s_1(z), s_2(z)$ are not identical zeros, then
$s(z)$ is not a unimodular constant. Therefore,
$D_{s,s(t_i)}(t_i)\ne 0,$ by the third assertion in Lemma
\ref{L:7.1}. Thus,
$$
0<D_{s,s(t_i)}(t_i)<\infty.
$$
If the second limit in (\ref{12.2100}) is also finite, then
$$
\lim_{z\to t_i}{\mathcal E} (z)=s(t_i)^*,
$$
by the first assertion in Lemma \ref{L:7.1}. Therefore,
(\ref{6.10}) turns into (\ref{7.12}) in this case. If
$D_{\cE,{s(t_i})^*}(t_i)$ is infinite, then, in view of
(\ref{12.2200}), the denominator in (\ref{6.10}) tends to
$\infty.$ Since the numerator $\cE(z)s(t_i)$ is bounded, the limit
in (\ref{6.10}) is $0.$ Thus, (\ref{7.12}) holds in this case also.
Theorem follows. \qed

\medskip

{\bf Proof of Statement 2 in Theorem \ref{T:1.5}:}
As it was already pointed out, $w$ is a solution of
Problem~\ref{P:8.1} if and only if it is of the form
(\ref{7.11}) with some (uniquely determined) parameter
$\cE\in\cS$ and satisfies
$$
\delta_{w,i}:=\gamma_i-d_{w,n_i}(t_i)= 0\quad (i=1,\ldots,k).
$$
The formula for $\delta_{w,i}$ is given in (\ref{7.12}) and it is
easily seen that $\delta_{w,i}=0$ if and only if
$$
D_{\cE, \, {s(t_i})^*}(t_i)+D_{s, \, s(t_i)}(t_i)=\infty.
$$
Since $D_{s, \, s(t_i)}(t_i)<\infty$ (by Corollary \ref{C:12.1}),
the latter is equivalent to $D_{\cE,\, {s(t_i})^*}(t_i)=\infty$
which happens, by the second assertion in Lemma \ref{L:7.1}, if
and only if either
$$
\liminf_{z\to t_i} \frac{1-|\cE(z)|^2}{1-|z|^2}=\infty ,
$$
or the function $\cE$ fails to have the nontangential limit ${
s(t_i})^*$ at $t_i$.\qed

\medskip

Note that vanishing of the gap at the point $t_i$ depends on
the local behavior of the parameter $\cE$ at this point only. The
number $s(t_i)^*$ absorbs all the interpolation data, though. Note
also that the maximum value of the gap $\delta_{w,i}$ is assumed
when $D_{\cE,\ {s(t_i})^*}(t_i)=0,$ what happens if and only if
$\cE(z)\equiv {s(t_i})^*.$

\bibliographystyle{plain}

\bigskip
MSC: 30E05, 446E22, 47A57

\end{document}